\documentclass[10pt]{article}

\usepackage{imsart}
\usepackage{amsmath,amssymb,amsfonts,amsthm,euscript,amsbsy,subfigure}
\RequirePackage[dvips,colorlinks]{hyperref}

\usepackage[dvips]{graphicx}
\graphicspath{{./}{figures/}}
\usepackage{amssymb,epsfig}
\usepackage{amsfonts}

\include{references}

\startlocaldefs
\numberwithin{equation}{section}
\theoremstyle{plain}

\endlocaldefs

\newtheorem{theorem}{Theorem}[section]
\newtheorem{proposition}[theorem]{Proposition}
\newtheorem{lemma}[theorem]{Lemma}

\newtheorem{remark}[theorem]{Remark}

\newcommand{\se}{\sigma}
\newcommand{\LL}{L}
\newcommand{\IE}{\mathbb{E}}

\newcommand{\Ebold}{\mathbb{E}}
\newcommand{\Pbold}{\mathbb{P}}

\newcommand{\Fcal} {{\cal F}}
\newcommand{\Ocal} {{\cal O}}
\newcommand{\nnb}       {\nonumber}

\begin{document}
\begin{frontmatter}
\title{Fixation Probability for Competing Selective Sweeps\protect}

\runtitle{Competing Selective Sweeps}

\begin{aug}
\author{\fnms{Feng} \snm{Yu}\ead[label=e3]{fyz@stats.ox.ac.uk}},
\author{\fnms{Alison} \snm{Etheridge}\ead[label=e2]{etheridg@stats.ox.ac.uk}},
and
\author{\fnms{Charles} \snm{Cuthbertson}\ead[label=e1]{cuthbert@stats.ox.ac.uk}
\ead[label=u1,url]{www.foo.com}}
\affiliation{University of Bristol, University of Oxford, and Morgan Stanley, UK
\footnote{CC supported by EPSRC DTA, FY supported by EPSRC/GR/T19537 while
at the University of Oxford}
}
\end{aug} 

\runauthor{Yu, Etheridge \& Cuthbertson}

\begin{abstract}
We consider a biological population in which a beneficial mutation is
undergoing a selective sweep when a second beneficial mutation arises at a
linked locus and we investigate the probability that both mutations will
eventually fix in the population.  Previous work has dealt with the case
where the second mutation to arise confers a smaller benefit than the first.
In that case population size plays almost no role. Here we consider the
opposite case and observe that, by contrast, the probability of both
mutations fixing can be heavily dependent on population size. Indeed the key
parameter is $\rho N$, the product of the population size and
the recombination rate
between the two selected loci.  If $\rho N$ is small, the probability that
both mutations fix can be reduced through interference to almost zero while
for large $\rho N$ the mutations barely influence one another.  The main
rigorous result is a method for calculating the fixation probability of a
double mutant in the large population limit.
\end{abstract}


\end{frontmatter}

\pagestyle{headings}
\section{Introduction}
\label{sec:intro}
Natural populations incorporate beneficial mutations through a combination
of chance and the action of natural selection.
The process whereby a beneficial mutation arises
(in what is generally assumed to be a large and otherwise neutral population)
and eventually spreads to the entire population is called a
{\em selective sweep}.
When beneficial mutations are rare, we can make the simplifying assumption that
selective sweeps do not overlap. A great deal is known about such
isolated selective sweeps (see e.g. Chapter 5 of Ewens~1979).\nocite{ewens1979}
Haldane~(1927) showed that under a discrete generation haploid model,
the probability that a beneficial allele \nocite{haldane:1927}
with selective advantage $\se$ eventually {\em fixes} in a population of
size $2N$,
i.e. its frequency increases from $1/(2N)$ to 1, is approximately $2\se$.
Much less is understood when selective sweeps overlap, i.e. when further
beneficial mutations arise at different loci during the timecourse of
a sweep.

Our aim here is to investigate the impact of the resulting interference in the
case when two sweeps overlap. In particular, we shall investigate the
probability that both beneficial mutations eventually become fixed
in the population. Because genes are organised on
chromosomes and chromosomes are in turn grouped into individuals,
different genetic loci do not evolve independently of one another. However,
in a dioecious population (in which chromosomes are carried in pairs), nor are
chromosomes passed down as intact units.  A given chromosome is inherited from
one of the two parents, but {\em recombination} or {\em crossover} events
can result in the allelic types at two distinct loci being inherited one from
each of the corresponding pair of chromosomes in the parent.
We refer to these chromosomes as `individuals'.

Each individual in the population will have a type denoted $ij$ where
$i,j\in \{0,1\}$. We use the first and second digit, respectively, to
indicate whether the individual carries the more recent or the older
beneficial mutation, and assume that the fitness effects of these two mutations
are additive. Suppose that a single advantageous allele
with selective advantage
$\se_1$ arises in an otherwise neutral (type 00) population of size $2N$,
corresponding to a diploid population of size $N$.
We use $X_{ij}$ to denote the proportion of individuals of type $ij$,
then the frequency of the favoured allele, $X_{01}$, will be
well-approximated by the solution to the stochastic differential equation
\begin{equation}
        dX_{01}=\se_1 X_{01}(1-X_{01}) \ ds
		+\sqrt{\frac{1}{2N}X_{01}(1-X_{01})} \ dW(s), \label{eq:X01s}
\end{equation}
where $s$ is the time variable,
$\{W(s)\}_{s\geq 0}$ is a standard Wiener process, and
$X_{01}(0)=1/(2N)$ (Ethier \& Kurtz~1986, Eq. 10.2.7).
\nocite{ethier/kurtz:1986}
If the favoured allele reaches frequency $p$, then the probability that it
ultimately fixes is $$\frac{1-e^{-2N\se_1p}}{1-e^{-2N\se_1}}.$$
If a sweep {\em does} take place then (conditioning on fixation) we obtain
\begin{equation*}
        d\tilde X_{01}=\se_1 \tilde X_{01}(1-\tilde X_{01})
			\coth(N\se_1\tilde X_{01}) \ ds
                +\sqrt{\frac{1}{2N}\tilde X_{01}(1-\tilde X_{01})} \ dW(s)
\end{equation*}
and from this it is easy to calculate the expected duration of the sweep.
Writing
$\tilde T_{fix}=\inf\{\left. s\geq 0: \tilde X_{01}(s)=1 \right|
                \tilde X_{01}(0)=1/(2N) \}$,
we have (see for example Etheridge et al.~2006)
\nocite{EPW06}
\begin{equation}
\label{eq:fix_time}
        \IE[\tilde T_{fix}]=\frac 2 {\se_1} \log (2N\se_1)
                +\mathcal{O}\left(\frac{1}{\se_1}\right)
\end{equation}
and the variance $var[\tilde T_{fix}]$ is $\mathcal{O}(1/\se_1^2)$.
More generally, an analogous Green function calculation to that leading to
equation~(\ref{eq:fix_time}) gives that the expected time for the selected
locus to reach frequency $\epsilon(N)$ is
$\log (2N\se_1\epsilon(N) )/\se_1+\mathcal{O}(1/\se_1)$.
This is the same as the expected time for $\tilde X_{01}$ to increase from
$1-\epsilon (N)$ to 1. On the other hand, for $\delta=\mathcal{O}(1)$,
the time for $\tilde X_{01}$ to increase from $\delta$ to $1-\delta$ is
$\mathcal{O}(1/\se_1)$.  As a result, for large populations, during almost
all the timecourse of the sweep $\tilde X_{01}$ is either close to zero
or close to one.

Now suppose that during the selective sweep of type 01 described
by~(\ref{eq:X01s}), more specifically, when $X_{01}$ reaches a level $U$,
another beneficial mutation with selection coefficient
$\se_2$ occurs at a second linked locus in a randomly chosen individual,
and the recombination rate between these two loci is $\rho$.
If we assume that the arrival time of the second mutation is
uniformly distributed over the timecourse of the sweep of the first mutation
and that $N$ is large, then we can expect either $U$ or $1-U$
to be close to 0 but $\gg 1/(2N)$.
The new mutation can arise in a type 00 or 01 individual, forming
a single type 10 individual in the former case,
and a 11 individual in the latter case.
If the second mutation arises during the first half (in terms of time)
of the sweep of the first mutation, then $U$ is likely to be very small and
it is more likely for a type 10 individual to be formed.
Otherwise, the second mutation arises during the second half of the sweep
and the formation of a type 11 individual is more likely.

The case of the second beneficial mutation forming a type 11 individual
is relatively straightforward. Since type 11 is fitter than all other types,
its fixation is almost certain once it becomes `established' in the population,
i.e. when the number of type 11 individuals is much larger than 1.
If the population size is very large, then it only takes a short time
to determine whether type 11 establishes itself, and we can assume the
proportion of type 01 individuals remains roughly constant during this time.
Hence the fixation probability of type 11 is essentially its
establishment probability, which is approximately $2(\se_2+\se_1(1-U))$,
twice the `effective' selective advantage of type 11 in a population
consisting of $2NU$ type 01 and $2N(1-U)$ type 00 individuals.

The case of the second beneficial mutation forming a type 10 individual
is far more interesting. In order for both mutations to sweep through
the population, recombination must produce an individual carrying
both mutations. The relative strength of selection acting on the two loci now
becomes important.
The case of $\se_1>\se_2$ has been dealt with in Barton~(1995)
and Otto \& Barton~(1997). \nocite{Barton95} \nocite{otto/barton:1997}
Here, since type 01 is already present in significant numbers when the new
mutation arises (and type 01 is fitter than type 10), the trajectory of
$X_{01}$ is well approximated by the logistic growth curve
$1/(1+\exp(-\se_1 t))$ until $X_{11}$ reaches a level of $\Ocal(1)$.
At that point, fixation of type 11 is all but certain.
Barton~(1995) then uses a branching process
approximation to estimate the establishment probability of a type 11
individual produced by recombination. In particular, his
approach is {\em independent} of population size. Not surprisingly,
he finds that the
fixation probability of the second mutation is reduced if it arises
as a type 10 individual, but increased if it arises as a type 11 individual.
Simulation studies performed in Otto \& Barton~(1997) confirm these
findings in the case $\se_1>\se_2$.

Gillespie~(2001) considers the effects of repeated substitutions at a
strongly selected locus on a completely linked (i.e. there is no recombination)
{\em weakly selected} locus, extending his work in Gillespie~(2000),
where he considers a linked {\em neutral} locus. He too sees little dependence
of his results on population size, leading him to suggest repeated genetic
hitchhiking events as an explanation for the apparent insensitivity of
the genetic diversity of a population to its size.
Kim~(2006) extends the work of Gillespie~(2001)
by considering the effect of repeated sweeps on a tightly (but
not completely) linked locus. This whole body of work is concerned, in our
terminology, with $\se_1>\se_2$.
\nocite{kim:2006}
\nocite{gillespie:2000}
\nocite{gillespie:2001}

The case of $\se_2>\se_1$ brings quite a different picture. The analysis
used in Barton~(1995) breaks down for the following reason: because the second
beneficial mutation is more competitive than the first, type 10 is destined
to start a sweep itself if it gets established in the population.
Once $X_{10}$ reaches $\Ocal(1)$, $X_{01}$ is no longer well approximated
by a logistic growth curve and in fact will decrease to 0.
The fixation probability of type 11 will then depend on the nonlinear
interaction of all four types, $\{11, 10, 01, 00\}$, and our analysis will
show that it is {\em heavily dependent} on population size.
See Figure~\ref{fig:time_avg} below.

\begin{figure}[h!]
\centering
\includegraphics[height = 2.3in, width = 2.8in]{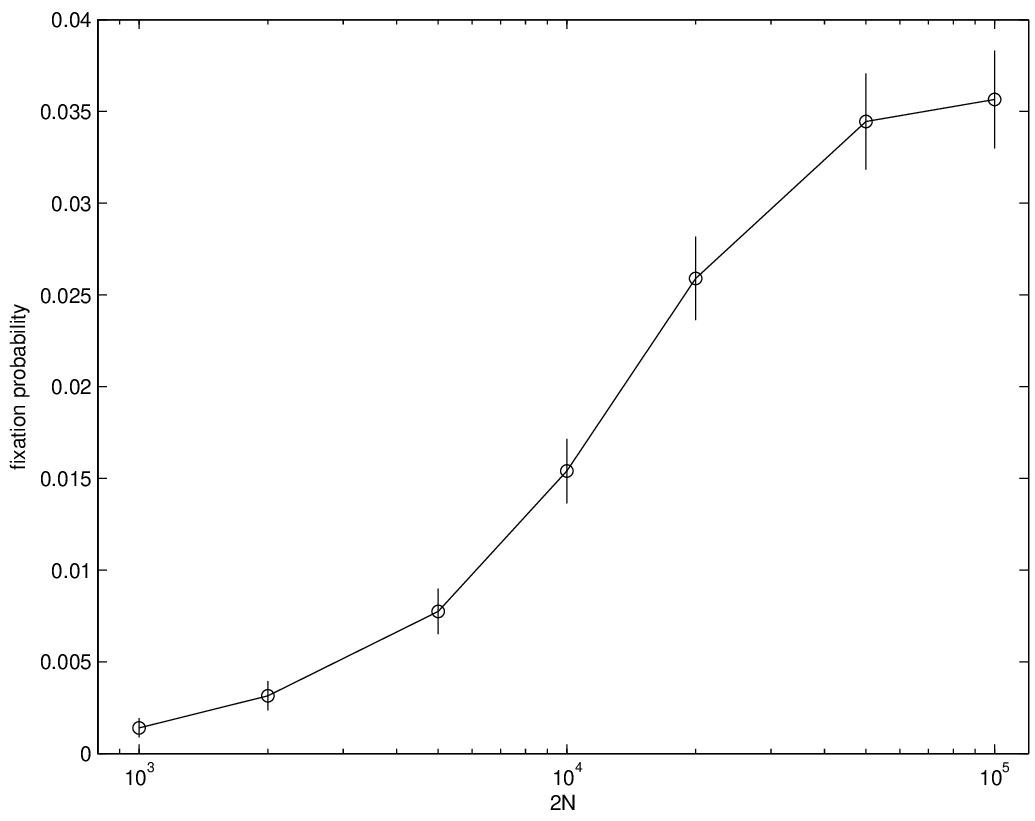}
\caption{Simulation results for fixation probability of type 11
for the following initial condition:
the second mutation arises in a type 00 individual, when $(2N)^{0.7}$
individuals in the population has the first mutation (i.e. are of type 01).
Vertical bars denote two standard deviations.
Parameter values: $\se_1=0.012$, $\se_2=0.02$, $\rho=4\times 10^{-5}$
(recombination coefficient).}
\label{fig:time_avg}
\end{figure}

This paper is organized as follows. 
In \S\ref{sec:model} we set up a continuous time Moran model for the 
evolution of our population.  In the biological literature, it would be
more usual to consider a Wright-Fisher model, in which the population evolves
in discrete, non-overlapping generations.
The choice of a Moran model, in which generations overlap,
is a matter of mathematical convenience. One expects similar results
for a Wright-Fisher model. The choice of a discrete individual based model
rather than a diffusion is forced upon us by our method of proof, but is
anyway natural in a setting where population size plays a r\^ole in the 
results. A brief analysis of our model, for very large $N$,
leads to our main rigorous result, Theorem~\ref{thm:main},
which provides a method to calculate the asymptotic ($N\to\infty$)
fixation probability of type 11 when $\se_2>\se_1$. We
discuss the case of moderate $N$ in \S\ref{sec:modN}.
The rest of the paper is devoted to proofs, with \S\ref{sec:proofs} containing
the proof of Theorem~\ref{thm:main} and \S\ref{sec:case1} containing
the proof of Proposition~\ref{prop:Z10}. Results in \S\ref{sec:case1}
rely on supporting lemmas of \S\ref{sec:lemmas}.

\section{Main Results}
\label{sec:mainresults}
\subsection{A Moran Model for Two Competing Selective Sweeps}
\label{sec:model}
In this section we describe our model for the evolution of two competing
selective sweeps.  We use the notation from the introduction for the four 
possible types of individual in the population
$I = \{00, 10, 01, 11\}$,
and assume that at the time when the second mutation arises, the number 
$U\in \{0,1,\ldots ,2N\}$ of type 01 individuals in the population is known.
From now on we use
$t=0$ to denote the time when the second mutation arises.
As explained in \S\ref{sec:intro}, we may assume that
$U$ is much larger than 1. 

Let $\sigma\in [0,1]$ be the selective advantage of the
second beneficial mutation and 
$\sigma\gamma$ be the selective advantage of the 
first beneficial mutation (for some $\gamma >0$).
The recombination rate between the two selected loci is denoted by
$\rho$ which we assume to be $o(1)$.
We use $\{(\eta_n \zeta_n),n=1,\ldots,2N\}$ to denote the types of individuals
in the population. At time $t=0$, we assume that the population of $2N$
individuals consists of $2N - U - 1$ type 00 individuals,
$U$ type 01 individuals and 1 type 10 individual.
The dynamics of the model are as follows:
\begin{enumerate}
\item[1.] {\em Recombination}: Each ordered pair of individuals,
	$(\eta_m \zeta_m)$ and $(\eta_n \zeta_n) \in I$,
	is chosen at rate $\rho/(2N)$.
	With probability $1/2$, $(\eta_m \zeta_n)$ replaces $(\eta_m \zeta_m)$.
	Otherwise, $(\eta_n \zeta_m)$ replaces $(\eta_m \zeta_m)$.  

\item[2.] {\em Resampling (and selection)}: Each ordered pair of individuals,
	$(\eta_m \zeta_m)$ and $(\eta_n \zeta_n) \in I$,
	is chosen at rate $1/(2N)$.
	With probability $p(\eta_m \zeta_m, \eta_n \zeta_n)$ given by
\begin{eqnarray}
	p(ij,kl) : = \frac{1}{2}(1 + \sigma (i - k) + \sigma \gamma(j - l)),
		\nnb
\end{eqnarray}
a type $(\eta_m \zeta_m)$ individual replaces $(\eta_n \zeta_n)$.
Otherwise a type $(\eta_n \zeta_n)$ individual replaces $(\eta_m \zeta_m)$.
\end{enumerate}

\begin{remark} Evidently we must assume $\sigma(1+\gamma)\le 1$ to ensure
that all probabilities used in the definition of the model are in $[0,1]$.
\end{remark}

\begin{remark}
If $\rho$ and $\sigma$ are small, then
decoupling recombination from the rest of the reproduction process
does not affect the behaviour of the model a great deal and it will simplify
analysis.
\end{remark}


Let $\Pbold$ denote the law of this Moran particle system, and
$r^+_{ij}$ and $r^-_{ij}$ be the rates at which $X_{ij}$
increases and decreases by $1/(2N)$, respectively, then
\begin{eqnarray*}
	r^+_{10} &=& N X_{10} [(1+\sigma)(1-X_{10})
                        -\sigma(1+\gamma) X_{11} - \sigma\gamma X_{01}] \\
                && \qquad \qquad \qquad
                + \rho N (2X_{11}X_{00}+X_{10}X_{11}+X_{10}X_{00}) \\
        r^-_{10} &=& N X_{10} [(1-\sigma)(1-X_{10})
                        +\sigma(1+\gamma) X_{11} + \sigma\gamma X_{01}] \\
                && \qquad \qquad \qquad
                + \rho N X_{10} (X_{00}+2X_{01}+X_{11})
\end{eqnarray*}
\vspace{-1.0cm}
\begin{eqnarray*}
	r^+_{01} &=& N X_{01} [(1+\sigma\gamma)(1-X_{01})
                - \sigma(1+\gamma)X_{11} - \sigma X_{10}] \\
        && \qquad \qquad \qquad
                + \rho N (X_{00}X_{01}+X_{11}X_{01}+2X_{11}X_{00})
                \\
        r^-_{01} &=& N X_{01} [(1-\sigma\gamma)(1-X_{01})
                + \sigma(1+\gamma)X_{11} + \sigma X_{10}] \\
        && \qquad \qquad \qquad
                + \rho N X_{01} (X_{00}+2X_{10}+X_{11})
\end{eqnarray*}
\vspace{-0.9cm}
\begin{eqnarray*}
	r^+_{11} &=& N X_{11} [(1+\sigma(1+\gamma))(1-X_{11})
                        -\sigma X_{10} - \sigma\gamma X_{01}] \nnb \\
        && \qquad \qquad \qquad
                + \rho N (2X_{10}X_{01}+X_{10}X_{11}+X_{01}X_{11})
                \nnb \\
        r^-_{11} &=& N X_{11} [(1-\sigma(1+\gamma))(1-X_{11})
                        +\sigma X_{10} + \sigma\gamma X_{01}] \nnb \\
        && \qquad \qquad \qquad
                + \rho N X_{11} (2X_{00}+X_{01}+X_{10}) \nnb
\end{eqnarray*}
\vspace{-1.0cm}
\begin{eqnarray}
	r^+_{00} &=& N X_{00} [1-X_{00}
                        -\sigma(1+\gamma) X_{11} -\sigma X_{10}
                        -\sigma\gamma X_{01}] \nnb \\
		&& \qquad \qquad \qquad
                + \rho N (X_{01}X_{00}+X_{00}X_{10}+2X_{01}X_{10}) \nnb \\
        r^-_{00} &=& N X_{00} [1-X_{00}
                        +\sigma(1+\gamma) X_{11} +\sigma X_{10}
                        + \sigma\gamma X_{01}] \nnb
                \\&& \qquad \qquad \qquad
                + \rho N X_{00} (X_{01}+2X_{11}+X_{10}).  \label{eq:rij}
\end{eqnarray}

\subsection{Analysis and Results for Large $N$}
\label{sec:results}
We are concerned primarily with the case of very large population sizes,
which is the regime where our main rigorous result,
Theorem~\ref{thm:main}, operates. A non-rigorous analysis for moderate
population sizes
based on very similar ideas is also possible but will appear
in Yu \& Etheridge~(2008).
\nocite{TwoSweepsNonrig}

To motivate our result, we present a heuristic analysis of the possible
scenarios. The proof of our main result fills in the necessary steps
to make this rigorous. If the second beneficial mutation gives rise
to a single type 10 individual, then the process whereby type 11 becomes fixed
must proceed in three stages and our
approach is to estimate the probability of each of these hurdles being 
overcome.  First, following the appearance of the new mutant, $X_{10}$ 
must `become established', by which we mean achieve appreciable frequency in
the population.  Without this, there will be no chance of step two: 
recombination of a type 01 and a type 10 individual to produce a type 11.
Finally, type 11 must become established (after which its ultimate
fixation is essentially certain).  Of course this may not happen the {\em first}
time a new recombinant is produced.  If type 11  becomes extinct and 
neither $X_{01}$ nor $X_{10}$ is one, then we can go back to step two.

We assume the first mutation has been undergoing a selective sweep prior
to the arrival of the second mutation. Before the arrival of the second
beneficial mutation (during which $X_{10}$ and $X_{11}$ are both 0),
we can write
\begin{eqnarray*}
        X_{01}(s) = \frac 1 {2N} + M_{01}(s)
        	+ \int_0^s \sigma\gamma X_{01}(u) (1-X_{01}(u)) \ du,
\end{eqnarray*}
where $M_{01}$ is a martingale with maximum jump size $1/(2N)$ and
quadratic variation $\langle M_{01} \rangle(s) = \frac{1+\rho}{2N}
        \int_0^s X_{01}(u)(1-X_{01}(u)) \ du$.
i.e.  $\left< M_{01} \right>$ is the unique previsible process
such that $M_{01}(s)^2-M_{01}(0)^2-\langle M_{01}\rangle(s)$
is a martingale. See e.g. \S~II.3.9 of Ikeda \& Watanabe~(1981).
\nocite{IkedaWatanabe}
We drop the martingale term $M_{01}$ and approximate the trajectory of
$X_{01}$ using a logistic growth curve, i.e.
$X_{01}(s)\approx 1/(1+(2N-1)\exp(-\sigma\gamma s))$ which solves
$\frac{d X_{01}}{ds} = \sigma\gamma X_{01}(s) (1-X_{01}(s))$
and $X_{01}(0)=1/(2N)$.
As discussed in \S\ref{sec:intro},
if we assume that the arrival time of the second mutation is uniformly
distributed on the timecourse of the sweep of the first
and $N$ is large, then $X_{01}$ spends most of the time near 0 or near 1.

We divide into two cases.
\begin{enumerate}
\item[1.] The second mutation arises during the first half of the sweep
	of the first mutation, i.e. when $X_{01}<1/2$.
\item[2.] The second mutation arises during the second half of the sweep
	of the first mutation, i.e. when $X_{01}\ge 1/2$.
\end{enumerate}

In Case 2, $X_{01}$ is close to 1 and it is most likely that
the second mutation arises in a type 01
individual to form a single type 11 individual, in which case
the fixation probability is roughly the same as the establishment
probability of type 11 arising in a population consisting entirely
of type 01 individuals, which in turn is roughly $2\sigma/(1+\sigma)$.

From now on, we focus on the more interesting Case 1. In what follows,
$t=0$ will be the time of arrival of the second beneficial mutation.
There it is most likely that the second mutation arises in a type 00
individual resulting in a single type 10 individual in the population.
If we approximate the growth of $X_{01}$ by a logistic growth curve, then
it reaches $1/2$ at time
$\frac 1 {\sigma\gamma} \log(2N-1) \approx \frac 1 {\sigma\gamma} \log(2N)$.
Choosing the time of the introduction of the new mutation uniformly on
$[0,\frac 1 {\sigma\gamma} \log(2N)]$ we see that at $t=0$,
$X_{01} \approx (2N)^{-\zeta}$, where $\zeta \sim Unif[0,1]$.

The establishment probability for type 10 in this case
is relatively easy to estimate. Since $\se_2 > \se_1$, type 10 either dies
out becomes established
before $X_{01}$ can grow to be a significant proportion
of the population. Therefore the establishment probability of type 10 is
almost the same as a type 10 arising in a population consisting
entirely of type 00 individuals, roughly $2\sigma/(1+\sigma)$.

We observe that if type 11 does get established, then since it has fitness
advantage over all other types, the probability that it eventually fixes
is very close 1 (this follows from Lemma~\ref{lem:estfix}).
Therefore we can concentrate on the behaviour of $X$ before
$X_{11}$ reaches say $(\log (2N))/(2N)$,
which is still very small compared to 1.
After type 10 is established and prior to type 11 being established,
we approximate $X_{10}$ and $X_{01}$ deterministically.
Until either $X_{10}$ or $X_{01}$ is $\Ocal(1)$, both grow
roughly exponentially, so assuming that type 10 gets established, we have
\begin{eqnarray}
	X_{10}(t) \approx \frac 1 {2N} e^{\sigma t}, \
	X_{01}(t) \approx \frac 1 {(2N)^\zeta} e^{\sigma\gamma t}.
	\label{approx:X1001}
\end{eqnarray}
We divide Case 1 further into two sub-cases. See Figure~\ref{fig:10} for
an illustration.
{
\begin{figure}[h!]
\centering
\subfigure[Case 1a: $\zeta = 0.3, \gamma=0.6$] {
\includegraphics[height = 2.0in, width = 2.0in]{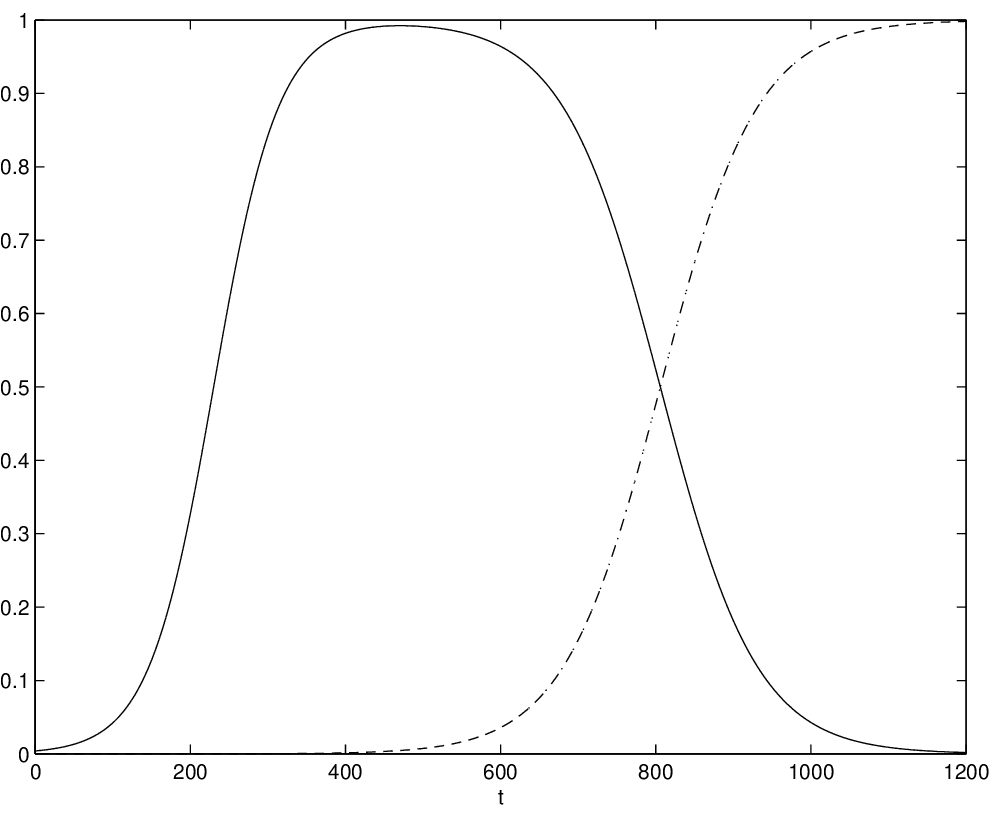}}
\subfigure[Case 1b: $\zeta = 0.7, \gamma=0.6$] {
\includegraphics[height = 2.0in, width = 2.0in]{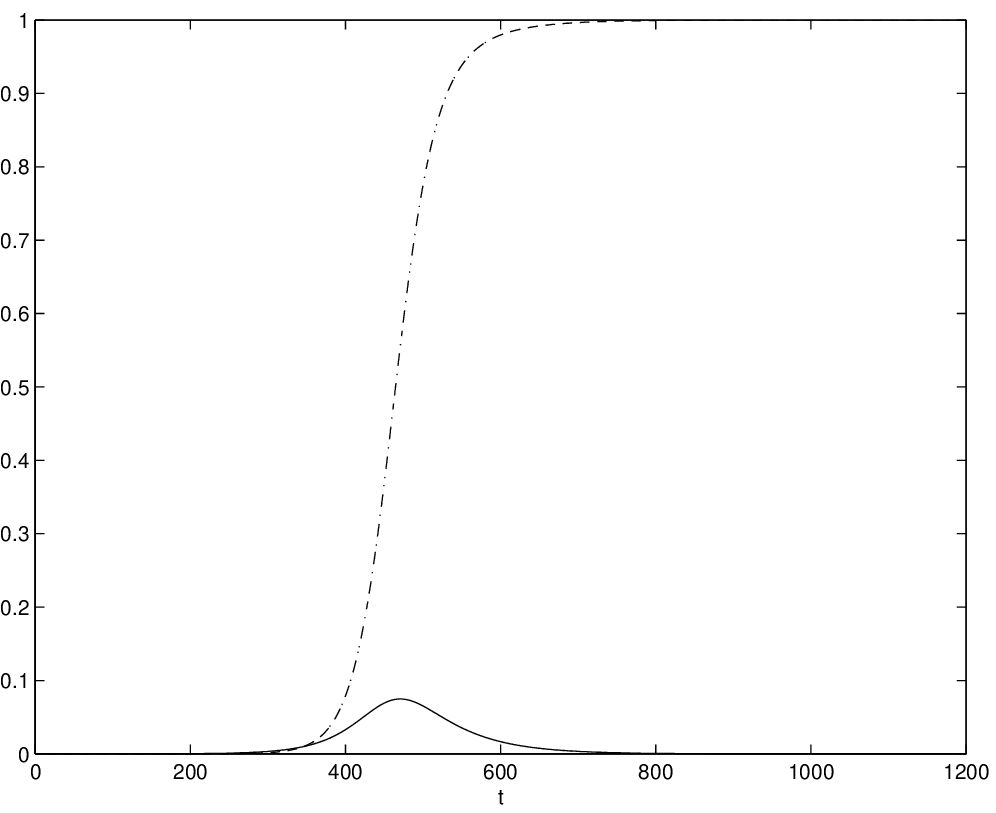}}
\caption{Approximate trajectories of $X_{01}$ (solid line)
and $X_{10}$ (dashed line) when $X_{11}$ is small: these curves
are obtained assuming they undergo deterministic logistic growth with initial
condition $X_{10}(0)=(2N)^{-1}$ and $X_{01}(0)=(2N)^{-\zeta}$.
 Parameter values: $\sigma=0.02$, $(2N)=10^8$. In
Case 1a, $X_{01}$ reaches almost 1 before being displaced by $X_{10}$,
but in Case 1b, $X_{01}$ never reaches $\Ocal(1)$.}
\label{fig:10}
\end{figure}
}

\noindent {\bf Case 1a, $\zeta<\gamma$.}
The approximation~(\ref{approx:X1001}) fails once either $X_{10}$ or $X_{01}$
reaches $\Ocal(1)$, which occurs at time
$\frac 1 {\sigma} \log (2N) \wedge \frac \zeta {\sigma\gamma} \log (2N)$.
If $\zeta<\gamma$, then $X_{01}$ reaches $\Ocal(1)$ before $X_{10}$, and
will further increase to almost 1 (which takes time only $\Ocal(1)$) before
$X_{10}$ reaches $\Ocal(1)$. At this time, which we denote $T_1$,
the population consists almost entirely of types 01 or 10.
Type 10, already established but
still just a small proportion of the population, will then proceed to
grow logistically, displacing type 01 individuals until $X_{10}$ is
close to 1 at time $T_2$. During $[T_1,T_2]$ (of length $\Ocal(1)$),
both $X_{01}$ and $X_{10}$ are $\Ocal(1)$,
so we expect $\Ocal(\rho N)$ recombination events between them producing
$\Ocal(\rho N)$ type 11 individuals.  Each type 11 individual has a
probability of at least $2\sigma\gamma/(1+\sigma\gamma)$
of eventually becoming the common ancestor of all individuals in the
population. So if we want to get a nontrivial limit (as $N\to\infty$)
for the fixation probability of type 11, we should take $\rho=\Ocal(1/N)$.
When we use the term nontrivial here, we mean that as $N\to\infty$, (i)
the fixation probability does not tend to 0, due to a lack of recombination
events between type 10 and type 01 individuals, and (ii) nor does it tend
to the establishment probability of type 10, due to infinitely
many type 11 births, one of which is bound to sweep to fixation.

\noindent {\bf Case 1b, $\zeta>\gamma$.}
In this case, $X_{10}$ reaches $\Ocal(1)$ at time roughly
$\frac 1 {\sigma} \log (2N)$, before $X_{01}$ does, and $X_{01}$ is
$\Ocal((2N)^{\gamma-\zeta})$ at this time. Furthermore, the biggest $X_{01}$
can get is $\Ocal((2N)^{\gamma-\zeta})$ since $X_{10}$ will very soon
afterwards increase to almost 1, after which $X_{01}$ will exponentially
decrease (since type 01 is less fit than type 10). Hence we expect
$\Ocal(\rho N^{1+\gamma-\zeta})$ recombination events between type 10
and type 01, and the `correct' scaling for $\rho$ is
$\rho=\Ocal(N^{\zeta-\gamma-1})$ in this case.

In case 1a, we take $\rho=\Ocal(1/N)$, then most of the recombination
events between type 10 and type 01 individuals occur when type 10
is logistically displacing type 01, i.e. in the time interval $[T_1,T_2]$.
During this time, we can approximate $X_{10}$ and $X_{01}$ by
$Z_{10}$ and $1-Z_{10}$, respectively, where $Z_{10}$ is {\em deterministic}
and obeys the logistical growth equation with parameter
$\sigma(1-\gamma)$, twice the advantage of type 10 over type 01.
We can further approximate $X_{11}$ by a birth and death process $Z_{11}$ with
deterministic but time-varying rates that depend on $Z_{10}$. Specifically,
the rates of increase and decrease for $Z_{11}$ are the same as
$r^\pm_{11}$ in~(\ref{eq:rij}), but with
$X_{10}$ replaced by $Z_{10}$, $X_{01}$ replaced by $1-Z_{10}$
and $X_{11}$ replaced by 0.

The probability that $X_{11}$ gets established, i.e. reaches
\[ \delta_{11}=\lceil\log (2N)\rceil/(2N), \]
is then approximated by the probability that the birth and death process
$Z_{11}$ reaches $\delta_{11}$. The latter can be found by solving
the forward equation for the process $Z_{11}$, which can be found
in~(\ref{eq:p11}). We define the fixation time of the
Moran particle system of \S\ref{sec:model}:
\begin{eqnarray*}
        T_{fix} = \inf\{t \ge 0: X_{ij}(t) = 1
		\mbox{ for some } ij \in I\}.
\end{eqnarray*}
We observe that the Markov chain $(X_{00},X_{01},X_{10})$ has
finitely many states and the recurrent states are
$R=\{(0,0,0),(0,0,1),(0,1,0),(1,0,0)\}$. Every other state is transient and
there is positive probability of reaching $R$ starting from any transient
state in finite time. Therefore
\[ T_{fix}<\infty \ a.s. \]
Our main result, Theorem~\ref{thm:main} below, concerns Case 1a,
which is the most likely scenario if $\gamma$ is close to 1.

\begin{theorem}
If $\zeta<\gamma<1$ and $\rho=\Ocal(1/N)$, then
there exists $\delta>0$, whose value depends on $\rho$, $\sigma$, $\gamma$,
and $\zeta$, such that
\[ \left| \Pbold \left(X_{11}(T_{fix}) = 1 \right)
        - \frac{2\sigma}{1+\sigma} p^{(11)}_{\delta_{11}}(T_\infty)
        \right| \le N^{-\delta}
\]
for sufficiently large $N$, where $p^{(11)}(t)$ solves
the forward equation~(\ref{eq:p11}).
\label{thm:main}
\end{theorem}

In the above, $\frac{2\sigma}{1+\sigma}$ corresponds to the
establishment probability of type 10, while $p^{(11)}_{\delta_{11}}(T_\infty)$
approximates the establishment probability of type 11 conditioning on
type 10 becoming established.
Figure~\ref{fig:11fix} compares fixation probabilities obtained from
simulation, our non-rigorous calculation (which we briefly discuss
in \S\ref{sec:modN} below), and the large population limit of
Theorem~\ref{thm:main}. In Figure~\ref{fig:11fix}(a)
we hold $\rho N$ constant in this simulation,
and observe that the fixation probability of type 11 increases but
does not change drastically as $N$ becomes large. The reason for the
drop in the fixation probability of type 11 when $N$ is small
may be because in this case, the early phase for $X_{01}$ is very short
and hence grows quickly to reduce the establishment probability of type 10.
In Figure~\ref{fig:11fix}(a), we use a population size of $2N=50,000$
to approach the large population limit of Theorem~\ref{thm:main}. At
$2N=50,000$, it takes roughly 12 hours on a PC to obtain one data point
in Figure~\ref{fig:11fix}, which is run with 20,000 realisations.
Apparently this population size still results in underestimates of the limiting
large population limit.
{
\begin{figure}[h!]
\centering
\subfigure[] {
\includegraphics[height = 2.3in, width = 2.3in]{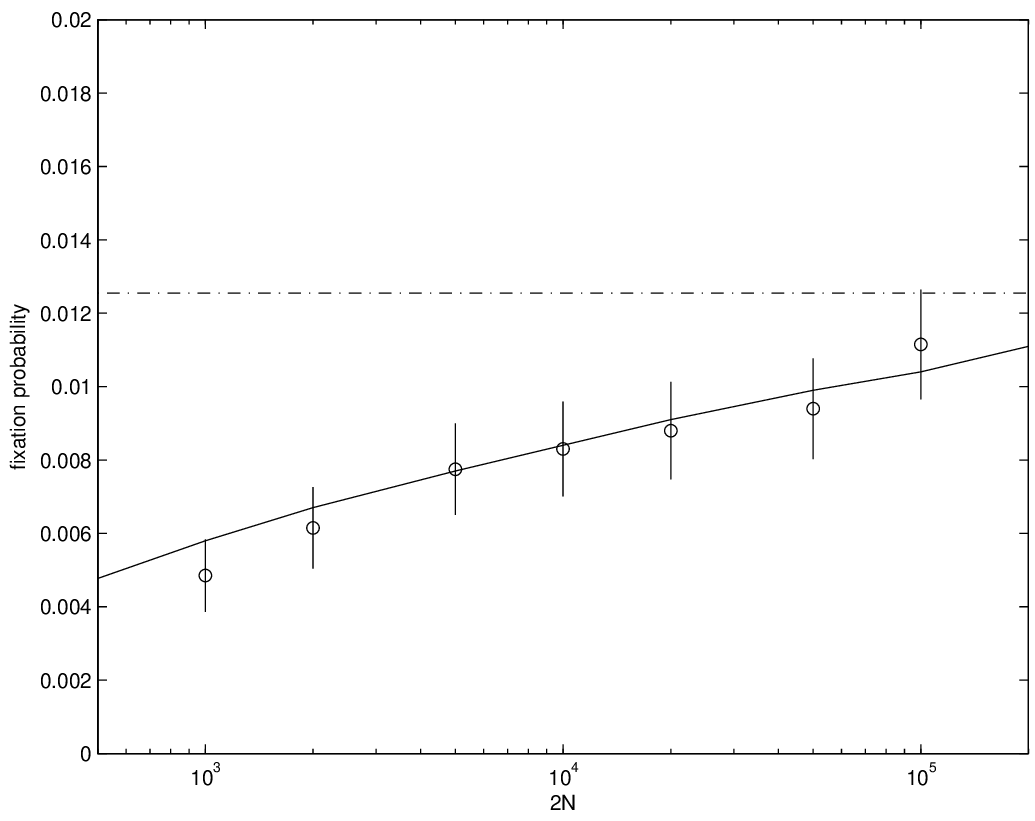}}
\subfigure[] {
\includegraphics[height = 2.3in, width = 2.3in]{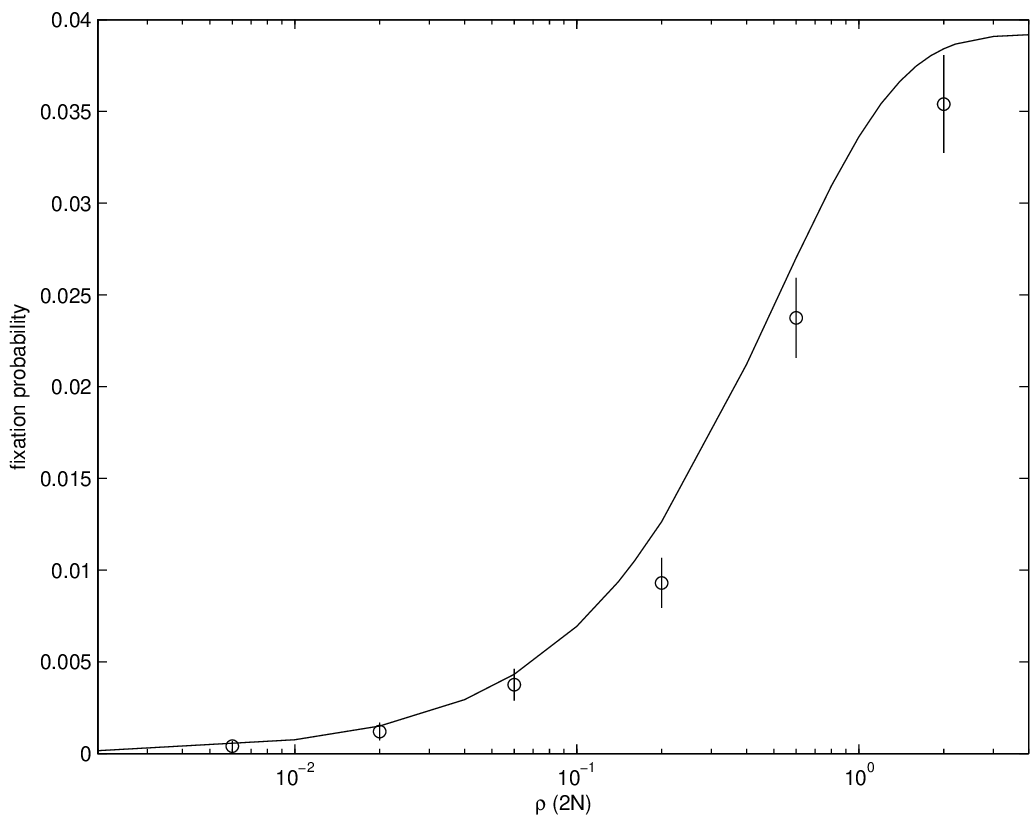}}
\caption{Fixation probability of type 11: circles denote data points
from simulations with vertical bars denoting one standard deviation.
(a) varying population size: 
the solid line denotes probabilities obtained using our non-rigorous
calculation, and the dashed line denotes the large population limit of
Theorem~\ref{thm:main}, with $\rho (2N)=0.2$.
(b) varying $\rho (2N)$: the solid line plots the large
population limit of Theorem~\ref{thm:main}, and the simulation
uses population size $2N=50,000$.
Other parameter values:
$\sigma=0.02$, $\zeta=0.3$
and $\gamma=0.6$. }
\label{fig:11fix}
\end{figure}
}

We expect a similar result for Case 1b, for which we provide an outline here.
We take $\epsilon\le (\gamma-\zeta)/(2+\gamma)$ and
$t_1=\frac{1-\epsilon}{\sigma} \log (2N)$, then at time $t_1$, we expect
$X_{10}$ to be either 0 (with probability approximately
$\frac{1-\sigma}{1+\sigma}$, as in Case 1a) or $\Ocal((2N)^{-\epsilon})$ and
$X_{01}$ to be roughly $(2N)^{(1-\epsilon)\gamma-\zeta}\le (2N)^{-2\epsilon}$.
Since $X_{01}$ and $X_{11}$ can be expected to be quite small before $t_1$,
they exert little influence on the trajectory of $X_{10}$, which jumps by
$\pm 1/(2N)$ at roughly the following rates:
\begin{eqnarray*}
        r^+_{10} \approx N (1+\sigma+\rho) X_{10}, \
        r^-_{10} \approx N (1-\sigma+\rho) X_{10}.
\end{eqnarray*}
Hence before $t_1$, $2N X_{10}$ resembles a continuous-time branching process
$Z$ with generating function of offspring distribution in the form of
$u(s)=\frac 1 2 (1+\sigma+\rho) s^2 + \frac 1 2 (1-\sigma+\rho) - (1+\rho) s$.
Using Theorem III.8.3 of Athreya \& Ney (1972), we can calculate
$E[e^{-u W}]$ for $W=\lim_{t\to\infty} e^{-\sigma t} Z(t)$ and conclude that
$W$ is distributed according to
$\frac{1-\sigma+\rho}{1+\sigma+\rho} \delta_{0}(x)
	+ \exp(\frac{2\sigma}{1+\sigma+\rho} x) \ dx$ for $x\ge 0$.
Hence the conditional distribution function of $X_{10}(t_1)| X_{10}(t_1)>0$
resembles $Exp(\frac{1+\sigma+\rho}{2\sigma} (2N)^{-\epsilon})$,
an exponential distribution with mean
$\frac{1+\sigma+\rho}{2\sigma} (2N)^{-\epsilon}$, as $N\to\infty$.

From time $t_1$ onwards, until either $X_{10}$ gets very close to 0 or $X_{01}$
becomes much smaller than $\Ocal((2N)^{(1-\epsilon)\gamma-\zeta})$,
we can assume that the paths of
$X_{01}$  and $X_{10}$ resembles those of $Z_{01}$ and $Z_{10}$, respectively,
where
\begin{eqnarray*}
	d Z_{10} &=& Z_{10} [(1+\sigma)(1-Z_{10}) - \sigma\gamma Z_{01}]
		\ dt \\
	d Z_{01} &=& Z_{01} [(1+\sigma\gamma)(1-Z_{01}) - \sigma Z_{10}]
		\ dt
\end{eqnarray*}
with the initial condition $Z_{10}(t_1)$ drawn according to 
$Exp(\frac{1+\sigma+\rho}{2\sigma} (2N)^{-\epsilon})$ and
$Z_{01}(t_1)=(2N)^{(1-\epsilon)\gamma-\zeta}$.
As in Case 1a, we can then approximate $X_{11}$ by a birth and death
process $Z_{11}$ with rates the same as $r^\pm_{11}$ from~(\ref{eq:rij})
but with $X_{10}$ replaced by $Z_{10}$ and $X_{01}$ replaced by $Z_{01}$.
The probability that $Z_{11}$ reaches $\delta_{11}$ can then be found
by solving the forward equation for $Z_{11}$.
Finally, we integrate this probability against all initial conditions
for $Z_{10}$, drawn according to
$Exp(\frac{1+\sigma+\rho}{2\sigma} (2N)^{-\epsilon})$.
The proof of such a result is more tedious than that of Theorem~\ref{thm:main}
but makes use of similar ideas.

\subsection{Brief Comment on Moderate $N$}
\label{sec:modN}
For moderate population sizes,
the observation in Case 1a of \S\ref{sec:results} that $X_{01}$ increases
to close to 1 before $X_{10}$ reaches $\Ocal(1)$ breaks down.
We can, however, compute the distribution function $f_T$ of the random time
$T_{10;\delta_{10}}$ when $X_{10}$ hits a certain level $\delta_{10}$,
assuming that $X_{01}$ grow logistically before $T_{10;\delta_{10}}$.
From $T_{10;\delta_{10}}$ onwards and before $X_{11}$ hits $\delta_{11}$,
$X_{10}$ grows roughly deterministically,
displacing both type 10 and type 00, so we can approximate $X_{11}$
by $Z_{11}$, a birth and death process with time-varying jump rates
in the form of $r^\pm_{11}$ in~(\ref{eq:rij}), but with
$X_{10}$, $X_{01}$ and $X_{00}$ replaced by their deterministic approximations.
Assuming $T_{10;\delta_{10}}=t$, we can numerically solve the forward
equation for $Z_{11}$, which is directly analogous to~(\ref{eq:p11}),
to find the probability that $Z_{11}$ eventually
hits $\delta_{11}$, which we denote by $p^{(11)}_{est}(t)$.
The dependence of $p^{(11)}_{est}$ on $t$ comes through the initial
condition $X_{01}$ for the ODE system, which depends on $T_{10;\delta_{10}}$.
The fixation probability of type 11 is then approximately
$\int p^{(11)}_{est}(t) f_T(t) \ dt$. This is the algorithm we use
to produce the solid line in Figure~\ref{fig:11fix}(a) and is given in its full
detail in Yu \& Etheridge~(2008).

\section{Proof of the Main Theorem}
\label{sec:proofs}
We first define some of the functions, events, and stochastic processes needed
for the proof, then give some intuition, before we proceed with
the proof of Theorem~\ref{thm:main}.
We begin by describing a deterministic process $Y_{10}$ and
a birth and death process $Y_{11}(t)$ which,
up to a shift by a random time, are $Z_{10}$ and $Z_{11}$ described in
\S\ref{sec:results}, respectively. They approximate the trajectories
of $X_{10}$ and $X_{11}$, respectively, after the establishment of type 10.
To describe the (time-inhomogeneous) rates we need the solution
\begin{eqnarray}
	\LL(t;y_0,\theta) = \left[1+\left(\frac 1 {y_0} -1\right)
	e^{-\theta t}\right]^{-1} \label{def:L}
\end{eqnarray}
to the logistic growth equation $\LL(t;y_0,\theta)
= y_0 + \theta \int_0^t \LL(s;y_0,\theta)(1-\LL(s;y_0,\theta)) \ ds$.
In what follows, $a_0=\zeta/(3\gamma)$ is a constant,
$c_1$, $c_2$, $c_3$ are constants (slightly smaller than $\Ocal(1)$)
that we specify precisely in Proposition~\ref{prop:Z10}, and
\begin{eqnarray}
	t_0 &=& \frac{a_0}{\sigma} \log (2N), \
        t_{early} = \frac{1.01\log (2N)}{\sigma(1-\gamma)-\rho},
		\label{def:t0} \\
	t_{mid} &=& \frac 1 {\sigma(1-\gamma)} \log \frac{1-c_1}{c_1}, \
	t_{late} = \frac{1.02}{\sigma\gamma} \log (2N). \nnb 
\end{eqnarray}
These deterministic times roughly correspond to the lengths of
the `stochastic', `early' (an upper bound), `middle', and `late' phases of
$X_{01}$, whose r\^ole is described in more detail in \S\ref{sec:case1}.
During the time interval when $Y_{10}$
is between $c_1$ and $1-c_1$, whose length is exactly $t_{mid}$,
there are birth events of $Z_{11}$ corresponding roughly
to recombination events between type 10 and 01 individuals.
For $t\in [0,t_{mid})$, we define
\begin{eqnarray*}
	Y_{10}(t) &=& \LL(t;c_1,\sigma(1-\gamma)) \\
        \beta^+(z,t) &=& N z [(1+\sigma(1+\gamma))(1-z)
                - (\sigma-\rho) Y_{10}(t)
                - (\sigma\gamma-\rho) (1-Y_{10}(t))] \\
        && \qquad + 2\rho N Y_{10}(t) (1-Y_{10}(t)) \\
        \beta^-(z,t) &=& N z [(1-\sigma(1+\gamma)+2\rho)(1-z)
                + (\sigma-\rho) Y_{10}(t) \\
        && \qquad + (\sigma\gamma-\rho) (1-Y_{10}(t))],
\end{eqnarray*}
and for $t\ge t_{mid}$, we define
\begin{eqnarray*}
        Y_{10}(t) &=& 1 \\
        \beta^+(z,t) &=& N (1+\sigma\gamma+\rho) z (1-z) \\
        \beta^-(z,t) &=& N (1-\sigma\gamma+\rho) z (1-z).
\end{eqnarray*}
We then take $Y_{11}$ to be a birth and death process with birth
and death rates $\beta^+ (Y_{11},t)$ and $\beta^- (Y_{11},t)$, respectively
(i.e. $Y_{11}$ jumps by $\pm 1/(2N)$ at rates
$\beta^+(Y_{11},t)$ and $\beta^+(Y_{11},t)$, respectively),
and initial condition $Y_{11}(0)=0$. It is absorbed on hitting $\delta_{11}$.

It is convenient to write $k_-=k-1/(2N)$ and $k_+=k+1/(2N)$.
$Y_{11}$ is run until time $t_{mid}+t_{late}$.
The probability that $Y_{11}$ hits $\delta_{11}$ before then
can be found by solving a system of ODE's.
Let $p^{(11)}$ satisfy
\begin{eqnarray*}
        \frac{d}{dt} p^{(11)}_k(t)
        = \beta^+ \left(k_-,t \right) p^{(11)}_{k_-}(t)
                + \beta^- \left(k_+,t \right) p^{(11)}_{k_+}(t)
        	- (\beta^+(k,t)+\beta^-(k,t)) p^{(11)}_k(t)
\end{eqnarray*}
for $k=1/(2N),\ldots,\delta_{11,-}$ where
$\delta_{11,-}=\delta_{11}-1/(2N)$, and
\begin{eqnarray}
        \frac{d}{dt} p^{(11)}_0(t)
        &=& \beta^- \left(1/(2N),t \right) p^{(11)}_{1/(2N)}(t)
                - \beta^+(0,t) p^{(11)}_0(t) \nnb \\
        \frac{d}{dt} p^{(11)}_{\delta_{11}}(t)
        &=& \beta^+ \left(\delta_{11,-},t \right)
                        p^{(11)}_{\delta_{11,-}}(t)
                - \beta^-(\delta_{11},t) p^{(11)}_{\delta_{11}}(t)
	\label{eq:p11}
\end{eqnarray}
with initial condition $p^{(11)}_k(0) = {\bf 1}_{\{ k = 0 \} }$. Then
\begin{eqnarray}
	\Pbold( Y_{11} \mbox{ hits } \delta_{11} \mbox{ before }
	t_{mid}+t_{late} ) = p^{(11)}_{\delta_{11}}(t_{mid}+t_{late}).
	\label{eq:p11Y}
\end{eqnarray}

We use the following convention for stopping times:
\begin{eqnarray}
        T_{ij;x} &=& \inf\{t\ge 0: X_{ij}\ge x \}, \
        T_{Z;x} =  \inf\{t\ge 0: Z\ge x \} \label{eq:Tij} \\
        S_{Y,Z,diff} &=& \inf\{ t\ge 0: Y(t)\ne Z(t) \} \nnb
\end{eqnarray}
for any $ij\in \{00,01,10,11\}$ and processes $Y$ and $Z$,
and define stopping times
\begin{eqnarray*}
        T_\infty &=& T_{10;c_1}+t_{mid}+t_{late}, \\
        S_{10,01,rec} &=& \inf\{ t\ge 0:
                \mbox{ there is a recombination event}\\
        && \mbox{between a type 10 and a type 01 individual
		before time $t$} \}.
\end{eqnarray*}
We define events
\begin{eqnarray*}
	E_1 &=& \{ X_{10}(t_0) > 0 \} \\
	E_2 &=& \{ T_{10;c_1} \le T_{11;1/(2N)} \wedge (t_0+t_{early}) \}
        	\cap \{ X_{01}(T_{10;c_1}) \ge 1-c_1-c_2 \} \\
        E_3 &=& \{ |X_{10}(t)-Z_{10}(t)| \le c_3
        	\mbox{ and } X_{00}(t) \le \sqrt{c_2} \mbox{ for all } \\
        && \qquad t \in [T_{10;c_1},
                        T_{Z_{10};1-c_1}\wedge T_{11;\delta_{11}}] \} \\
	E_4 &=& \{ T_{Z_{10};1-c_1} \le T_{11;\delta_{11}} \} \\
	E_5 &=& \{ X_{11}(t)+X_{10}(t) > 1-\sqrt{c_1}
                \mbox{ for all } t \ge T_{Z_{10};1-c_1} \} \\
	E_6 &=& \left\{ X_{11}(T_\infty)+X_{10}(T_\infty) = 1 \right\} \\
	E_7 &=& \{ X_{11}(t)= Z_{11}(t) \mbox{ for all }
                t\in [T_{10;c_1},T_\infty \wedge T_{11;\delta_{11}}] \} \\
        E_8 &=& \{ T_{11;\delta_{11}} \le T_\infty \mbox{ or }
                        X_{11}(T_\infty)=Z_{11}(T_\infty) = 0 \}.
\end{eqnarray*}
We observe that $T_{11;1/(2N)}\ge S_{10,01,rec}$. First we outline the
intuition behind these definitions:
$t_0$ is the length of the initial `stochastic' phase for $X_{10}$. At $t_0$,
with high probability $X_{10}$ either is $\Ocal((2N)^{a_0-1})$
or has hit 0 (event $E_1^c$). In the latter case, there is no need to
approximate $X_{10}$ any further. On the other hand, if $E_1$ occurs,
then type 10 is very likely to be established by $t_0$ and,
with high probability, grows almost deterministically to reach level $c_1$
(slightly smaller than $\Ocal(1)$) at time $T_{10;c_1}$.
Furthermore, as discussed in \S\ref{sec:intro}, in Case 1a,
since $\zeta<\gamma$, with high probability $X_{01}(T_{10;c_1})$
is close to 1. Hence conditional on $E_1$, the event $E_2$ is very likely.

For paths in $E_2 \cap E_1$, we define
\begin{eqnarray}
	Z_{10}(T_{10;c_1}+t) = Y_{10}(t), \ Z_{11}(T_{10;c_1}+t) = Y_{11}(t)
	\label{eq:ZY}
\end{eqnarray}
to be the approximations for the trajectories of $X_{10}$ and $X_{11}$,
respectively, from time $T_{10;c_1}$ onwards.
For convenience, we define $Z_{10}(t)=Z_{11}(t)=0$ for $t\le T_{10;c_1}$.
With the convention of~(\ref{eq:Tij}),
\[ T_{Z_{10};1-c_1} = T_{10;c_1} + t_{mid}, \]
and we observe that $Z_{10}(t)=1$ for $t\ge T_{Z_{10};1-c_1}$.
Since $X_{01}(T_{10;c_1})\approx 1$, $X_{00}(T_{10;c_1})$ is very small
and is unlikely to recover because type 00 is the least fit type.
During $[T_{10;c_1},T_{Z_{10};1-c_1}]$, with high probability,
type 10 grows logistically at rate $\sigma(1-\gamma)$, displacing
type 01. Hence conditional on $E_1\cap E_2$, $E_3$ is very likely.
During $[T_{10;c_1},T_{Z_{10};1-c_1}]$, the definition of $Z_{11}$
takes into account recombination
events between type 01 and 10 individuals that produce type 11
individuals at a rate of $\rho (2N) X_{01} X_{10}$, which in the definition
of $Z_{11}$, is approximated by $\rho (2N) Z_{10} (1-Z_{10})$. Notice
that we can approximate $X_{01}$ by $1-Z_{10}$ since we assume throughout that
$X_{11}\le\delta_{11}$, which is very small.
Outside the time interval $[T_{10;c_1},T_{Z_{10};1-c_1}]$, either $X_{10}$ is
very small or very close to 1 (which means $X_{01}$ is very small), hence
we ignore any recombination events. Because $Z_{11}$ closely approximates 
$X_{11}$, conditional on $E_3\cap E_2\cap E_1$,
event $E_7$ has a high probability.

After $T_{Z_{10};1-c_1}$, $X_{11}+X_{10}$ is likely to remain close to 1
(event $E_5$) and hit 1 at time $T_\infty$ (event $E_6$). We ignore
any more recombination events between type 10 and 01 and $Z_{11}$
is a time-changed branching process during this time.
If $Z_{11}$ has not hit $\delta_{11}$ by time $T_{Z_{10};1-c_1}$ (event $E_4$),
then we continue to keep track of $Z_{11}$ until $T_\infty$, at which
time it most likely has already hit either $\delta_{11}$ or 0 (event $E_8$).
In the latter case, we regard type 11 as having failed to establish and
since $X_{10}$ is most likely to be 1 (event $E_6$) at $T_\infty$,
the earlier mutation has gone extinct.
If $X_{11}$ hits $\delta_{11}$ before $T_\infty$, we regard
type 11 as having established and hence it will, with high probability,
eventually sweep to fixation (Lemma~\ref{lem:estfix}).

Proposition~\ref{prop:Z10} below estimates the probabilities
of events $E_1$ through $E_8$. These are `good' events, on which we
can approximate the establishment probability of type 11 by the probability
that $Z_{11}$ hits $\delta_{11}$ by time $T_\infty$.
Proposition~\ref{prop:Z10} is essential for the proof
of Theorem~\ref{thm:main}, and will be proved in \S\ref{sec:case1}.

\begin{proposition}
There exists positive constants $\delta_{10,3}$ and $\delta_{10,4}>0$
whose exact value depends on
$\sigma$, $\gamma$ and $\zeta$, such that $c_1, c_2, c_3$ in the definition
of $E_1,\ldots,E_8$ are all $\le N^{-\delta_{10,3}}$
and for sufficiently large $N$,
\begin{eqnarray*}
&(a)&	\left| \Pbold(E_1^c) - \frac{1-\sigma+\rho}{1+\sigma+\rho} \right|
        \le C_{\rho,\gamma,\sigma} N^{-\delta_{10,3}} \\
&(b)&	\Pbold (E_2^c \cap E_1) \le C_{\rho,\gamma,\sigma} N^{-\delta_{10,3}}
		\\
&(c)&	\Pbold(E_3^c \cap E_2 \cap E_1)
		\le C_{\rho,\gamma,\sigma} N^{-\delta_{10,3}} \\
&(d)&	\Pbold(E_5^c \cap E_4\cap E_3 \cap E_2 \cap E_1)
		\le C_{\rho,\gamma,\sigma} N^{-\delta_{10,3}} \\
&(e)&	\Pbold\left(E_6^c \cap E_4\cap E_3 \cap E_2 \cap E_1
		\right) \le C_{\rho,\gamma,\sigma} N^{-\delta_{10,3}}.
\end{eqnarray*}
Consequently, we have
$(f) \ \Pbold\left(E_6^c \cap E_4 \cap E_2 \cap E_1
                \right) \le C_{\rho,\gamma,\sigma} N^{-\delta_{10,3}}$.
Furthermore,
\begin{eqnarray*}
	&(g)& \Pbold (E_7^c \cap E_2\cap E_1 ) \le C_{\rho,\gamma,\sigma}
		(N^{-\delta_{10,3}}+N^{-\delta_{10,4}}) \\
        &(h)& \Pbold (E_8^c \cap E_7 \cap E_2\cap E_1 )
		\le C_{\rho,\gamma,\sigma} N^{-\delta_{10,4}}.
\end{eqnarray*}
\label{prop:Z10}
\end{proposition}

\begin{lemma}
$|\Pbold(X_{11}(T_{fix})\neq 1)-\Pbold(T_{11;\delta_{11}}<\infty)|
\le N^{\log \frac{1-\sigma\gamma+2\rho}{1+\sigma\gamma}}$.
\label{lem:estfix}
\end{lemma}
\proof
On $\{T_{11;\delta_{11}}<\infty \}$, $X_{11}$ dominates $\check X_{11}$,
a birth and death process with initial condition
$\check X_{11}(T_{11;\delta_{11}})=\delta_{11}=\lceil\log (2N)\rceil/(2N)$,
jump size $1/(2N)$, and the following jump rates
\begin{eqnarray*}
        \check r^+_{11} = N (1+\sigma\gamma)
		\check X_{11} (1-\check X_{11}), \
        \check r^-_{11} = N (1-\sigma\gamma+2\rho)
		\check X_{11} (1-\check X_{11}).
\end{eqnarray*}
Using standard Markov chain techniques, we may conclude
\[ \Pbold(\{ T_{\check X_{11};1} > T_{\check X_{11};0},
		T_{11;\delta_{11}}<\infty \})
	\le (2N)^{\log \frac{1-\sigma\gamma+2\rho}{1+\sigma\gamma}},
\]
which implies $\Pbold(\{X_{11}(T_{fix})\neq 1, \ T_{11;\delta_{11}}<\infty \})
\le (2N)^{\log \frac{1-\sigma\gamma+2\rho}{1+\sigma\gamma}}$.
Since $\{ X_{11}(T_{fix})= 1, \ T_{11;\delta_{11}}=\infty \}$ is
a set with probability 0, we have the desired result.
\qed

\vspace{.2cm}
\noindent
{\em Proof of Theorem~\ref{thm:main}.}
Recall from~(\ref{def:t0}) that $a_0=\zeta/(3\gamma)$
and $t_0 = \frac{a_0}{\sigma} \log (2N)$.
We first show that we can safely ignore $E_1^c$. Let
\begin{eqnarray*}
	E_9 &=& \left\{ X_{11}(t)=0 \mbox{ for all } t\le t_0 \right\}.
\end{eqnarray*}
Comparing with~(\ref{eq:rij}), we see that the jump process
$\hat X_{10}$ with initial condition $\hat X_{10}(0)=1/(2N)$,
jump size $1/(2N)$, and the following jump rates
\begin{eqnarray*}
        \hat r^+_{10} = N (1+\sigma) \hat X_{10}
		+ 3\rho N, \
        \hat r^-_{10} = N (1-\sigma) \hat X_{10}
\end{eqnarray*}
dominates $X_{10}$ for all time. Then
\begin{eqnarray*}
        d\hat X_{10} = dM + (\sigma \hat X_{10} + 1.5\rho) \ dt
\end{eqnarray*}
where $M$ is a martingale with maximum jump size $1/(2N)$ and quadratic
variation $\langle M \rangle$ satisfying
$d\langle M \rangle = \frac 1 {2N} (2\hat X_{10}+3\rho) \ dt$.
Hence
\begin{eqnarray*}
        E[\hat X_{10}(t)] = \left(\frac 1 {2N} + \frac{3\rho}{2\sigma} \right)
		e^{\sigma t} - \frac{3\rho}{2\sigma}
        \le \left(\frac 1 {2N} + \frac{3\rho}{2\sigma} \right) e^{\sigma t},
\end{eqnarray*}
We recall Burkholder's inequality in the following form:
\[ E\left[ \sup_{s\le t} |M(s)|^p \right] \le C_p \Ebold \left[
	\langle M \rangle(t)^{p/2} + \sup_{s\le t} |M(s)-M(s-)|^p \right],
\]
which may be derived from its discrete time version, Theorem 21.1 of
Burkholder~(1973). \nocite{Burkholder73}
We use this and Jensen's inequality to obtain
\begin{eqnarray}
        E\left[ \sup_{s\le t_0} |M(s)| \right]
	&\le&  E\left[ \sup_{s\le t_0} |M(s)|^2 \right]^{1/2} 
	\le \frac C N \left( 1 + N \int_0^{t_0}
			E[\hat X_{10}(s)+1.5\rho] \ ds \right)^{1/2} \nnb \\
        &\le& \frac C N + \frac{C_\sigma}{\sqrt N} \left( \rho t_0
                + (N^{-1} + \rho) e^{\sigma t_0} \right)^{1/2}
	\le C_{\rho,\sigma} N^{(a_0/2)-1}. \label{ineq:JB}
\end{eqnarray}
Therefore
\[ E\left[ \sup_{s\le t_0} \hat X_{10}(s) \right]
	\le E\left[ \sup_{s\le t_0} |M(s)| \right] + 1.5\rho t_0
		+ \sigma \int_0^{t_0} E[\hat X_{10}(s)] \ ds
	\le C_{\rho,\sigma} N^{a_0-1}.
\]
Since $\hat X_{10}$ dominates $X_{10}$, we have
\begin{eqnarray*}
        \Pbold \left(\sup_{s\le t_0} X_{10}(s)
                \ge (2N)^{2a_0-1}\right)
	\le C_{\rho,\sigma} N^{-a_0}.
\end{eqnarray*}
On $\{\sup_{s\le t_0} X_{10}(s) < (2N)^{2a_0-1} \}$,
the number of recombination events between type 10 and 01
during $[0,t_0]$ is at most $Poisson(2\rho (2N)^{2a_0-1} t_0)$,
hence
\[ \Pbold(E_9^c\cap E_1) \le \Pbold(E_9^c)
	\le C_{\rho,\sigma} (N^{-a_0} + N^{(2a_0-1)/2}) \]
for sufficiently large $N$.
On $E_9 \cap E_1^c$, type 10 has gone extinct by time $t_0$, before
a single individual of type 11 has been born, hence type 11 will not get
established, let alone fix. Therefore
\begin{eqnarray}
	\Pbold \left( \{ T_{11;\delta_{11}} <\infty \} \cap E_1^c \right)
	\le \Pbold (E_9^c \cap E_1^c)
	\le C_{\rho,\sigma} (N^{-a_0} + N^{(2a_0-1)/2}).
	\label{ineq:E8}
\end{eqnarray}

Now we concentrate on $E_1$ where type 10 has most likely established itself
at time $t_0$. The nontrivial event here is
$E_8 \cap E_7 \cap E_2 \cap E_1$.
Let $E_{81} = \{ T_{11;\delta_{11}} \le T_\infty \}$ and
$E_{82} = \{ T_{11;\delta_{11}} > T_\infty,
X_{11}(T_\infty)=Z_{11}(T_\infty) = 0 \}$, then
$E_8=E_{81}\cup E_{82}$. The following events have small probabilities
\begin{eqnarray}
	\Pbold(E_2^c \cap E_1)
	&\le& C_{\rho,\gamma,\sigma} N^{-\delta_{10,3}} \nnb \\
	\Pbold((E_8^c \cup E_7^c) \cap E_2 \cap E_1)
        &\le& C_{\rho,\gamma,\sigma} (N^{-\delta_{10,3}} + N^{-\delta_{10,4}})
		\nnb \\
	\Pbold(E_{82} \cap E_7 \cap E_6^c \cap E_2 \cap E_1)
	&\le& C_{\rho,\gamma,\sigma} N^{-\delta_{10,3}}, \label{ineq:E72}
\end{eqnarray}
by Prop~\ref{prop:Z10}(b), Prop~\ref{prop:Z10}(g-h), and
Prop~\ref{prop:Z10}(f), respectively, where the last estimate above
comes from the fact $E_{82}\subset E_4$.
There are two events with significant probabilities:
on $E_{82} \cap E_7 \cap E_6 \cap E_2 \cap E_1$, we have
$X_{11}(T_\infty) = 0, X_{10}(T_\infty) = 1$ hence
type 10 fixes by time $T_\infty$, and on
$E_{81} \cap E_7 \cap E_2 \cap E_1$,
$X_{11}=Z_{11}$ hits $\delta_{11}$ and get established by time $T_\infty$.
On both these events,
$X_{11}=Z_{11}$ until at least $T_\infty \wedge T_{11;\delta_{11}}$.
The union of these two events,
$E_{82} \cap E_7 \cap E_6 \cap E_2 \cap E_1$ and
$E_{81} \cap E_7 \cap E_2 \cap E_1$, and the three events in~(\ref{ineq:E72})
is $E_1$. On $E_1 \cap E_2$, for exactly one of the two events
$\{ T_{11;\delta_{11}} <\infty \}$
and $\{ T_{Z_{11};\delta_{11}} \le T_\infty \}$ to occur (i.e. either
the former occurs but the latter does not, or the latter occurs and the former
does not), one of the following three scenarios must occur:
\begin{enumerate}
\item $X_{11}$ and $Z_{11}$ disagree before $T_\infty$, i.e. $E_7^c$;
\item $X_{11}$ and $Z_{11}$ agree up to $T_\infty$, but do not hit
	$\{0,\delta_{11} \}$ before $T_\infty$, i.e. $E_8^c$;
\item $X_{11}$ and $Z_{11}$ agree up to $T_\infty$ and $X_{11}(T_\infty)=0$,
	but $X_{10}(T_\infty)<1$ thus allowing the possibility
	of type 11 being born due to recombination between type 10 and 01
	individuals after $T_\infty$, i.e. $E_6^c$.
\end{enumerate}
Hence
\begin{eqnarray*}
	\lefteqn{ \left| \Pbold \left( \{ T_{11;\delta_{11}} <\infty \}
			\cap E_1 \right)
		- \Pbold \left( \{ T_{Z_{11};\delta_{11}} \le T_\infty \}
			\cap E_1\right) \right| } \\
	&\le& \Pbold(E_2^c \cap E_1)
		+ \Pbold((E_8^c \cup E_7^c) \cap E_2 \cap E_1)
		+ \Pbold(E_{82} \cap E_7 \cap E_6^c \cap E_2 \cap E_1) \\
	&\le& C_{\rho,\gamma,\sigma} (N^{-\delta_{10,3}} + N^{-\delta_{10,4}})
\end{eqnarray*}
by~(\ref{ineq:E72}). From~(\ref{ineq:E8}), we have
\begin{eqnarray*}
	\lefteqn{ |\Pbold \left( T_{11;\delta_{11}} <\infty \right)
		- \Pbold \left( \{ T_{11;\delta_{11}} <\infty \} \cap E_1
		\right)| } \\
	&=& \Pbold \left( \{ T_{11;\delta_{11}} <\infty \}
			\cap E_1^c \right)
	\le C_{\rho,\sigma} N^{-a_0} + N^{(2a_0-1)/2}
\end{eqnarray*}
But by Proposition~\ref{prop:Z10}(a),
\[ \left| \Pbold(E_1) - \frac{2\sigma}{1+\sigma} \right|
        \le C_{\rho,\gamma,\sigma} N^{-\delta_{10,3}}. \]
We combine the three inequalities above to conclude
\begin{eqnarray*}
	\lefteqn{ \left| \Pbold \left( T_{11;\delta_{11}} <\infty \right)
	- \frac{2\sigma}{1+\sigma} \Pbold \left(
		T_{Z_{11};\delta_{11}} \le T_\infty | E_1\right) \right| } \\
	&\le& \left| \Pbold \left( T_{11;\delta_{11}} <\infty \right)
        	- \Pbold \left( T_{Z_{11};\delta_{11}} \le T_\infty
			| E_1\right) \Pbold(E_1) \right|
		+ C_{\rho,\gamma,\sigma} N^{-\delta_{10,3}} \\
	&=& \left| \Pbold \left( T_{11;\delta_{11}} <\infty \right)
                - \Pbold \left( \{ T_{Z_{11};\delta_{11}} \le T_\infty \}
                        \cap E_1\right) \right|
                + C_{\rho,\gamma,\sigma} N^{-\delta_{10,3}} \\
	&\le& N^{-\delta}
\end{eqnarray*}
for some $\delta>0$, and then use Lemma~\ref{lem:estfix},
as well as~(\ref{eq:p11Y}) and~(\ref{eq:ZY})  to obtain the desired conclusion.
\qed

\section{Proof of Proposition 3.1}
\label{sec:case1}
We divide the evolution of $X_{10}$ and $X_{01}$ roughly into
4 phases, `stochastic', `early', `middle', and `late', and use
Lemmas~\ref{lem:1dim_early},~\ref{lem:1dim_mid}, and~\ref{lem:1dim_late}
for each of the last 3 phases, respectively.
Lemma~\ref{lem:01early} deals with the early, middle, and late phases of
$X_{01}$. Because 
$X_{01}$ starts at $U = (2N)^{-\zeta} \gg 1/(2N)$ at $t=0$, it has no stochastic
phase. Its early phase is between $t=0$ and the time when $X_{01}$
reaches $c_{01,1}$. Its middle phase is between $c_{01,1}$ and
$1-c_{01,2}$, after which it enters the late phase.

For type 10, since $X_{10}(0)=1/(2N)$, whether it establishes itself is
genuinely stochastic (i.e. its probability tends to a positive constant
strictly less than 1 as $N\to\infty$).
The stochastic phase lasts for time $t_0$, when, with high probability,
either type 10 has established or it has gone extinct.
If $X_{10}$ reaches $\Ocal((2N)^{a_0-1})$ by time $t_0$, it enters 
the early phase, which is dealt with by Lemma~\ref{lem:10early}. Part (b)
of that lemma says that if $\zeta<\gamma$
(as mentioned before, we only deal Case 1a of \S\ref{sec:intro})
then it does not reach $c_{10,2}$ until
$X_{01}$ has entered its late phase, while part (c)
says that it does reach $c_{10,3}$ at some finite time. The proof of
Proposition~\ref{prop:Z10}(a-b)
reconciles various stopping times used in Lemmas~\ref{lem:01early}
and~\ref{lem:10early}, and prepares for part (c) of Proposition~\ref{prop:Z10},
which deals with the middle phase of $X_{10}$ during which $X_{10}$ increases
from $c_{10,3}$ to $1-c_{10,3}$, displacing $X_{01}$ in the process.
The $c_{ij,k}$'s we use throughout the rest of this paper are
small positive constants, all of $\Ocal((2N)^{-b_{ij,k}})$, whose exact values
are specified immediately below~(\ref{eq:constants}).

Recall the definition of the logistic growth curve $\LL(t;y_0,\theta)$
from~(\ref{def:L}). Throughout the rest of this section,
We use $\LL(t;(2N)^{-\zeta},\sigma\gamma)$ to approximate the trajectory of
$X_{10}$ during its early phase and
$t_{01;x}$ to denote the time when this approximation hits $x$, e.g.
$t_{01;c_{01,1}}$ below is when it hits $c_{01,1}$.
Furthermore, we use $t_{01,x,y}$ to denote the time this approximation
spends between $x$ and $y$. Thus
\begin{eqnarray*}
	\LL(t_{01;x};(2N)^{-\zeta},\sigma\gamma) = x \mbox{ and }
	\LL(t_{01,x,y};x,\sigma\gamma) = y.
\end{eqnarray*}
We also define
\[ t'_{01;1-c_{01,2}}=t_{01;c_{01,1}}+t_{01,0.9 c_{01,1},1-c_{01,2}}. \]
In the above, $t_{01,0.9 c_{01,1},1-c_{01,2}}$ is the length of time
for which we use the event $A_2$ in Lemma~\ref{lem:01early} below.
On the event $A_1$ defined in that lemma, $X_{01}$ reaches $0.9 c_{01,1}$
at time $t_{01;c_{01,1}}$, after which event $A_2$ ensures $X_{01}$ grows
to levels slightly smaller than $1-c_{01,2}$ after another time period
of length $t_{01,0.9 c_{01,1},1-c_{01,2}}$. Roughly speaking, the
time when $\LL(\cdot;(2N)^{-\zeta},\sigma\gamma)$ is between
$0.9 c_{01,1}$ and $c_{01,1}$ is counted twice.
We observe that
\begin{eqnarray}
	t_{01;c_{01,1}} &=& \frac 1 {\sigma\gamma} \log
		\frac{(2N)^\zeta-1}{\frac 1 {c_{01,1}} - 1} \nnb \\
	t'_{01;1-c_{01,2}} &=& t_{01;1-c_{01,2}} + t_{01,0.9 c_{01,1},c_{01,1}}
		\label{def:t01p} \\
	&=& \frac 1 {\sigma\gamma} \left\{ \log \left[ ((2N)^\zeta-1)
                \left(\frac 1 {c_{01,2}} - 1 \right) \right]
		+ \log \frac{\frac 1 {0.9 c_{01,1}}-1}{\frac 1 {c_{01,1}} - 1}
		\right\}. \nnb
\end{eqnarray}

We recall that $a_0 = \frac \zeta {3\gamma}$ and
define the constants required for the rest of the proof, as well
as $c_1$, $c_2$, and $c_3$ as required by Proposition~\ref{prop:Z10}:
\begin{eqnarray}
	a_1 &=& \frac \zeta {4\gamma} \wedge \frac{1-\zeta/\gamma} 4, \nnb \\
	b_{10,0} &=& a_0+a_1-1, \
	b_{10,2}=\frac{1-\zeta/\gamma} 2, \
	b_{10,3}=\frac{\gamma b_{10,2}}{90}, \nnb \\
	b_{01,0} &=& \frac \zeta 3, \
	b_{01,1}=b_{01,2}=\frac{\gamma b_{10,2}} 3, \nnb \\
	\delta_{01,1}  &=& \frac{\gamma b_{10,2}}{9}
		\le \frac \gamma 3 \left(b_{10,2} - b_{01,1}-b_{01,2}\right), \
        \delta_{10,2} = \frac{\delta_{01,1}}{60}=\frac{\gamma b_{10,2}}{540},
		\nnb \\
        \delta_{10,0}  &=&  2N c_{10,0} (c_{10,0}+c_{01,0}), \
        \delta_{10,1} = (a_0-a_1)/4, \nnb \\
	c_1 &=& c_{10,3}, \ c_2=(2N)^{-\delta_{01,1}/2}, \ c_3=(2N)^{-\delta_{10,2}}
	\label{eq:constants}
\end{eqnarray}
and $c_{ij,k}=(2N)^{-b_{ij,k}}$.
These choices imply $a_1 + b_{10,2} + b_{01,2}/\gamma < 1-\zeta/\gamma$,
which in turn implies the following:
\begin{eqnarray}
        (1-a_1) \log (2N) + \log c_{10,2} + \frac 1 \gamma \log c_{01,2}
		> \frac 1 \gamma \log ((2N)^\zeta-1), \nnb \\
	\log ((2N)^{1-a_1}-(2N)^{a_0})
                - \log \left(\frac{1}{0.9 c_{10,2}}-1\right)
                - \frac 1 \gamma \log \left(\frac 1 {c_{01,2}} - 1\right)
	\nnb \\
        > \frac 1 \gamma \log ((2N)^\zeta-1)
                + \frac 1 \gamma \log \frac 1 {0.9}, \nnb \\
        \log ((2N)^{1-a_1}-(2N)^{a_0})
                - \log \left(\frac{1}{0.9 c_{10,2}}-1\right)
        \qquad\qquad\qquad\qquad\qquad     \nnb \\
        \ge \frac 1 \gamma \left\{\log \left[ ((2N)^\zeta-1)
                \left(\frac 1 {c_{01,2}} - 1\right) \right]
                + \log \frac{\frac 1 {0.9 c_{01,1}}-1}{\frac 1 {c_{01,1}} - 1}
                \right\} = \sigma t'_{01;1-c_{01,2}}
	\label{ineq:constantsb}
\end{eqnarray}
for sufficiently large $N$. This will be needed in Lemma~\ref{lem:10early}.

\begin{lemma}
Let $R_{01}=T_{11;1/(2N)} \wedge T_{10;c_{10,2}}$. We define
\begin{eqnarray*}
	A_1 &=& \{X_{01}(s) \le 0.9 \LL(s;(2N)^{-\zeta},\sigma\gamma)
		\mbox{ for some } s\le t_{01;c_{01,1}} \wedge R_{01} \} \\
        A_2 &=& \{X_{01}(s) < \LL(s-t_{01;c_{01,1}};0.9 c_{01,1},\sigma\gamma)
                + (2N)^{-\delta_{01,1}} \mbox{ for some } \\
        && \qquad \qquad \qquad	s \in [t_{01;c_{01,1}},
                        t'_{01;1-c_{01,2}}\wedge R_{01}] \} \\
        A_3 &=& \{X_{10}(s)+X_{01}(s) \le 1-(2N)^{-\delta_{01,1}/2}
                \mbox{ for some } s\in [t'_{01;1-c_{01,2}},S_{10,01,rec}) \}.
\end{eqnarray*}
Then
\begin{eqnarray*}
&(a)& \Pbold(A_1) \le C_{\rho,\gamma,\sigma} N^{-(1-\zeta)/4} \\
&(b)&	\Pbold(A_2\cap A_1^c\cap \{t_{01;c_{01,1}}\le R_{01} \})
		\le (2N)^{-\delta_{01,1}} \\
&(c)&	\Pbold(A_3 \cap A_2^c \cap A_1^c
	\cap \{ t'_{01;1-c_{01,2}} \le R_{01} \}) \le C N^{-1/2}.
\end{eqnarray*}
Consequently,
\[ \Pbold((A_3 \cup A_2 \cup A_1)\cap \{ t'_{01;1-c_{01,2}} \le R_{01} \})
	\le C_{\rho,\gamma,\sigma} (2N)^{-\delta_{01,1}}. \]
\label{lem:01early}
\end{lemma}
\proof
\noindent {\bf Early Phase.}
Before the stopping time $R_{01}$,
the jump rates of $X_{01}$ satisfies
\begin{eqnarray*}
        r^+_{01} &\ge& N X_{01} [(1+\sigma\gamma+\rho)(1-X_{01})
                        - 1.1 \sigma c_{10,2}], \\
        r^-_{01} &\le& N X_{01} [(1-\sigma\gamma+\rho)(1-X_{01})
                        + 1.1 \sigma c_{10,2}].
\end{eqnarray*}
We take $\hat\xi=X_{01}$, $\alpha=1+\rho$, $\theta=\sigma\gamma$,
$\delta_0= 1.1 \sigma c_{10,2}$,
$\delta_1=c_{01,1}$, $\delta_2 = (1-\zeta)/4$, $Y$ such that
$Y(t)=(2N)^{-\zeta} + \int_0^t Y(s) (\sigma\gamma (1-Y(s))
	- 1.1 \sigma c_{10,2}) \ ds$,
and $u_0=\inf\{t: Y(t)=\delta_1\} > t_{01;c_{01,1}}$
in Lemma~\ref{lem:1dim_early} to obtain
\begin{eqnarray*}
        \Pbold \left( X_{01}(s) < 0.99 Y(s) \mbox{ for some }
		s\le t_{01;c_{01,1}} \wedge R_{01} \right)
        \le C_{\rho,\gamma,\sigma} N^{-(1-\zeta)/4}.
\end{eqnarray*}
Prior to $u_0$, $Y$ is sandwiched between
$\LL(\cdot;(2N)^{-\zeta},\sigma\gamma-1.2\sigma c_{10,2})$
and $\LL(\cdot;(2N)^{-\zeta},\sigma\gamma)$.
Since $\LL(t;(2N)^{-\zeta},\sigma\gamma)-\LL(t;(2N)^{-\zeta},\sigma\gamma-v)
\le (1-e^{-v t}) \LL(t;(2N)^{-\zeta},\sigma\gamma)$ for $v\le \sigma\gamma$,
we have
\begin{eqnarray*}
	Y(t) &\ge& \LL(t;(2N)^{-\zeta},\sigma\gamma-1.2\sigma c_{10,2})
	\ge e^{-1.2\sigma c_{10,2} t} \LL(t;(2N)^{-\zeta},\sigma\gamma) \\
	&\ge& 0.99 \LL(t;(2N)^{-\zeta},\sigma\gamma)
\end{eqnarray*}
for $t=\Ocal(\log N)$. Hence (a) follows.

\noindent {\bf Middle Phase.}
Before $R_{01}$, $X_{11}=0$.
Using the jump rates of $X_{01}$ in~(\ref{eq:rij}), we can write
\begin{eqnarray*}
        \lefteqn{ X_{01}(t\wedge R_{01}) = b_0 + M_{01}(t\wedge R_{01}) } \\
        && + \int_{u_1}^{t\wedge R_{01}} X_{01}(s)
                [\sigma\gamma (1-X_{01}(s)) - (\sigma+\rho) X_{10}(s)] \ ds,
\end{eqnarray*} 
where $M_{01}(\cdot\wedge R_{01})$ is a martingale with maximum jump size
$1/(2N)$ and quadratic variation
$\langle M_{01} \rangle(t\wedge R_{01}) = \frac{1+\rho}{2N}
	\int_{u_1}^{t\wedge R_{01}} X_{01}(s)(1-X_{01}(s)) \ ds$.
We apply Lemma~\ref{lem:1dim_mid} with $b_0=X_{01}(t_{01;c_{01,1}})$,
$u_1=t_{01;c_{01,1}}$, $u_2=t'_{01;1-c_{01,2}}$,
$\delta_1=b_{10,2}$, $\delta_2=\infty$,
$\epsilon_0=b_{01,1}$, $\epsilon_1=b_{01,2}$,
$\epsilon_3(t)=\epsilon_4(t)=0$, $T=R_{01}$ and $D_1=A_1^c$.
Then since $b_{10,2} > (b_{01,1}+b_{01,2}) \wedge \frac 1 2$,
we have
\begin{eqnarray*}
	\lefteqn{ \Pbold ( \{ |X_{01}(s) -
		\LL(s-t_{01;c_{01,1}};X_{01}(t_{01;c_{01,1}}),\sigma\gamma))|
		> (2N)^{-\delta_{01,1}} \mbox{ for some} } \\
        && s \in [t_{01;c_{01,1}}, t'_{01;1-c_{01,2}} \wedge R_{01}] \}
	\cap A_1^c\cap \{t_{01;c_{01,1}}\le R_{01} \} )
        \le (2N)^{-\delta_{01,1}},
\end{eqnarray*}
where $\delta_{01,1}$ is defined in~(\ref{eq:constants}).
Now for paths in $A_1^c\cap \{t_{01;c_{01,1}}\le R_{01} \}$, we have
$X_{01}(t_{01;c_{01,1}})\ge 0.9 c_{01,1}$ and hence
\[ \LL(s-t_{01;c_{01,1}};X_{01}(t_{01;c_{01,1}}),\sigma\gamma)
	\ge \LL(s-t_{01;c_{01,1}};0.9 c_{01,1},\sigma\gamma). \]
The desired conclusion in (b) follows.

\noindent {\bf Late Phase.}
On $A_1^c\cap A_2^c \cap \{ t'_{01;1-c_{01,2}} \le R_{01} \}$, 
since $\delta_{01,1}\le b_{01,2}$, we have
\begin{eqnarray*}
	X_{01}(t'_{01;1-c_{01,2}})
	&>& \LL(t'_{01;1-c_{01,2}}-t_{01;c_{01,1}};0.9c_{01,1},\sigma\gamma)
		-(2N)^{-\delta_{01,1}} \\
	&>& 1-c_{01,2}-(2N)^{-\delta_{01,1}}.
\end{eqnarray*}
Therefore $X_{00}(t'_{01;1-c_{01,2}})<2(2N)^{-\delta_{01,1}}$.
Before $S_{10,01,rec}$, $X_{11}=0$, and the jump rates of $X_{00}$ satisfy
\begin{eqnarray*}
        r^+_{00} \le N (1-\sigma\gamma+\rho) X_{00} (1-X_{00}), \
        r^-_{00} \ge N (1+\sigma\gamma+\rho) X_{00} (1-X_{00}).
\end{eqnarray*}
By Lemma~\ref{lem:1dim_late},
$\Pbold( \{ \sup_{t\ge t'_{01;1-c_{01,2}}} X_{00}(t) \le N^{-\delta_{01,1}/2}\}
\cap A_1^c\cap A_2^c \cap \{ t'_{01;1-c_{01,2}} \le R_{01} \}) \le C N^{-1/2}$,
which implies the desired conclusion in (c).
\qed

\vspace{.2cm}
For the remainder of this section, we define the following events
\begin{eqnarray*}
        A_{41} &=& \{ X_{10}(s) \ge (2N)^{a_0+a_1-1} \mbox{ for some }
                s \le t_0 \wedge T_{01;c_{01,0}} \wedge T_{11;1/(2N)} \} \\
        A_{42} &=& \{ X_{10}(t_0) \in [1,(2N)^{a_0-a_1-1}] \} \\
        A_4 &=& A_{41} \cup A_{42} \cup E_1^c \\
        B_4 &=&  \{t_0 \le T_{01;c_{01,0}} \wedge T_{11;1/(2N)} \} \\
        A_{51} &=& \{ X_{10}(s) \ge c_{10,2} \mbox{ for some }
                s\in [t_0,t'_{01;1-c_{01,2}} \wedge T_{11;1/(2N)}] \} \\
        A_{52} &=& \{ T_{10;c_{10,3}} \wedge T_{11;1/(2N)}
                        \ge t_0+t_{early} \}.
\end{eqnarray*}

\begin{lemma}
Recall that $t_0 = \frac{a_0}{\sigma} \log (2N)$,
$\delta_{10,0} = 2N c_{10,0} (c_{10,0}+c_{01,0})$,
$\delta_{10,1} = (a_0-a_1)/4$. We have
\begin{eqnarray*}
&(a)& \Pbold (A_{41}) \le 2 \delta_{10,0} t_0
	+ C_{\rho,\gamma,\sigma} N^{-a_1} \\
&& \Pbold( A_{42} \cap A_{41}^c \cap B_4)
        \le C_{\rho,\gamma,\sigma} N^{-a_1} + 2 \delta_{10,0} t_0 \\
&& \left| \Pbold(E_1^c \cap B_4)
                - \frac{1-\sigma+\rho}{1+\sigma+\rho} \right|
        \le 6\delta_{10,0} t_0 + C_{\rho,\gamma,\sigma} N^{-a_1}
                + \Pbold(B_4^c) \\
&(b)& \Pbold\left( A_{51} \cap A_4^c \cap B_4 \right)
	\le C_{\rho,\gamma,\sigma} N^{-\delta_{10,1}} \\
&(c)& \Pbold (A_{52} \cap A_4^c \cap B_4)
        \le C_{\rho,\gamma,\sigma} N^{-\delta_{10,1}}.
\end{eqnarray*}
\label{lem:10early}
\end{lemma}
\proof
\noindent {\bf Stochastic Phase.}
We define $R_{01,11}=T_{01;c_{01,0}} \wedge T_{11;1/(2N)}$.
Before $T_{11;1/(2N)}$, the jump rates of $X_{10}$ are
as follows:
\begin{eqnarray*}
        r^+_{10} &=& N X_{10} [(1+\sigma+\rho)(1-X_{10})
                        - (\sigma\gamma+\rho) X_{01}], \\
        r^-_{10} &=& N X_{10} [(1-\sigma+\rho)(1-X_{10})
                        + (\sigma\gamma+\rho) X_{01}].
\end{eqnarray*}
We define $\eta$ to be a jump process with $\eta(0)=1/(2N)$,
jump size $1/(2N)$ and jump rates as follows:
\begin{eqnarray*}
        r^+_{\eta,10} = N \eta (1+\sigma+\rho), \
        r^-_{\eta,10} = N \eta (1-\sigma+\rho).
\end{eqnarray*}
then prior to
$S_{X_{10},\eta,diff} \wedge T_{10;c_{10,0}} \wedge R_{01,11}$, we have
$|r^+_{10}-r^+_{\eta,10}| \le \delta_{10,0}$ and
$|r^-_{10}-r^-_{\eta,10}| \le \delta_{10,0}$.
Therefore $|X_{10}-\eta|$ is a jump process with initial value 0,
jump size $1/(2N)$ and jump rates at most $2\delta_{10,0}$,
and we can estimate the probability of $|X_{10}-\eta|$ becoming
nonzero before $t_0$:
\begin{eqnarray}
	\Pbold\left( S_{X_{10},\eta,diff} < t_0
        \wedge T_{10;c_{10,0}} \wedge R_{01,11}
        \right) \le 2 \delta_{10,0} t_0. \label{ineq:X10Y10}
\end{eqnarray}
Since $\eta$ is a branching process, Lemma~\ref{lem:branching}(a) implies
\begin{eqnarray*}
	\Pbold\left(\sup_{s\le t_0} \eta(s)\ge (2N)^{a_0+a_1-1}
		\right) &\le& C_{\rho,\gamma,\sigma} N^{-a_1} \\
	\Pbold(1\le \eta(t_0) \le (2N)^{a_0-a_1-1})
        &\le& C_{\rho,\gamma,\sigma} N^{-a_1} \nnb \\
	\left|\Pbold(\eta(t_0)=0)
		- \frac{1-\sigma+\rho}{1+\sigma+\rho}\right|
	&\le& \frac{1-\sigma+\rho}{1+\sigma+\rho} e^{-2\sigma t_0}
	\le (2N)^{-a_0}. \nnb
\end{eqnarray*}
Using~(\ref{ineq:X10Y10}), we can replace $\eta$ in the
above three estimates by $X_{10}$ if we allow an additional error term.
In particular,
\begin{eqnarray*}
	\lefteqn{ \Pbold\left(
		\sup_{s\le t_0\wedge R_{01,11}}
		X_{10}(s)\ge c_{10,0} = (2N)^{a_0+a_1-1} \right) } \\
	&=& \Pbold\left( \sup_{s\le t_0\wedge R_{01,11}}
                	X_{10}(s)\ge c_{10,0}, \
		S_{X_{10},\eta,diff} < t_0 \wedge T_{10;c_{10,0}}
			\wedge R_{01,11}  \right) \\
	&& + \Pbold\left(
                \sup_{s\le t_0\wedge R_{01,11}}
                        X_{10}(s)\ge c_{10,0}, \
                S_{X_{10},\eta,diff} \ge t_0 \wedge T_{10;c_{10,0}}
                        \wedge R_{01,11}  \right) \\
	&\le& 2 \delta_{10,0} t_0 + \Pbold\left(
                \sup_{s\le t_0\wedge R_{01,11}}
                        \eta(s)\ge c_{10,0} \right)
	\le 2 \delta_{10,0} t_0 + C_{\rho,\gamma,\sigma} N^{-a_1}.
\end{eqnarray*}
Similarly, we can obtain the second statement of (a) and
\begin{eqnarray*}
	\left| \Pbold(\{X_{10}(t_0)=0 \} \cap A_{41}^c \cap B_4)
		- \frac{1-\sigma+\rho}{1+\sigma+\rho} \right|
	\le 2 \delta_{10,0} t_0 + \Pbold(A_{41}) + \Pbold( B_4^c),
\end{eqnarray*}
which implies the third statement in (a).

\noindent {\bf Early Phase (Upper Bound).}
Before $T_{11;1/(2N)}$, the jump rates of $X_{10}$ satisfy
\[ r^+_{10} \le N (1+\sigma+\rho) X_{10} (1-X_{10}), \
        r^-_{10} \ge N (1-\sigma+\rho) X_{10} (1-X_{10}). \]
We take $\check\xi=X_{10}$, $\alpha=1+\rho$, $\theta=\sigma$,
$\delta_0=0$, $\delta_1=0.9 c_{10,2}$,
$\delta_2 = \delta_{10,1} = (a_0-a_1)/4$,
$Y(t) = \LL(t;X_{10}(t_0),2\sigma)$, and
$u_0=R_{10}=\inf\{t\ge 0: \LL(t;X_{10}(t_0),2\sigma) \ge \delta_1\}$
in Lemma~\ref{lem:1dim_early} to obtain
\begin{eqnarray*}
	\Pbold( \{ X_{10}(t_0+s) \ge 1.01 \LL(s;X_{10}(t_0),\sigma)
		\mbox{ for some } s \ge (R_{10} \wedge T_{11;1/(2N)})-t_0] \} \\
		\cap A_4^c \cap B_4 )
        \le C_{\rho,\gamma,\sigma} N^{-(1-\zeta)/4}.
\end{eqnarray*}
On $A_4^c\cap B_4 \subset
\{ X_{10}(t_0) \in ((2N)^{a_0-a_1-1},(2N)^{a_0+a_1-1}) \}$,
we have
\begin{eqnarray*}
        t_0+R_{10} &\ge& \frac 1 {\sigma} \left[ a_0 \log (2N)
                + \log ((2N)^{1-a_0-a_1}-1)
                - \log \left(\frac{1}{0.9 c_{10,2}}-1\right) \right] \\
        &\ge& t'_{01;1-c_{01,2}}
\end{eqnarray*}
by~(\ref{ineq:constantsb}) and the definition of $t'_{01;1-c_{01,2}}$
in~(\ref{def:t01p}). Hence if
\[ X_{10}(t_0)\in ( (2N)^{a_0-a_1-1}, (2N)^{a_0+a_1-1} ),\]
then $\LL(t'_{01;1-c_{01,2}}-t_0;X_{10}(t_0),\sigma)
\le \LL(R_{10};X_{10}(t_0),\sigma)
= 0.9 c_{10,2}$, which implies (b).

\noindent {\bf Early Phase (Lower Bound).}
Before $T_{11;1/(2N)}$, the jump rates of $X_{10}$ satisfy
\[ r^+_{10} \ge N (1+\sigma(1-\gamma)) X_{10} (1-X_{10}), \
        r^-_{10}\le N (1-\sigma(1-\gamma)+2\rho) X_{10} (1-X_{10}).\]
We take $\hat\xi$ to be $X_{10}$ shifted forward in time by $t_0$,
$\alpha=1+\rho$, $\theta=\sigma(1-\gamma)-\rho$,
$\delta_0=0$, $\delta_1=1.01 c_{10,3}$,
$\delta_2 = \delta_{10,1} = (a_0-a_1)/4$,
$Y(t)=\LL(t; N^{a_0-a_1-1},\sigma(1-\gamma)-\rho)$,
and $u_0=\inf\{t: Y(t)=1.01 c_{10,3} \}$
in Lemma~\ref{lem:1dim_early} to obtain
\begin{eqnarray*}
	\lefteqn{ \Pbold( \{ X_{10}(t_0+s) < 1.005
		\LL(s; N^{a_0-a_1-1},\sigma(1-\gamma)-\rho) } \\
 	&&\mbox{ for some } s \le u_0 \wedge (T_{11;1/(2N)}-t_0) \}
	\cap A_4^c \cap B_4)
        \le C_{\rho,\gamma,\sigma} N^{-\delta_{10,1}}.
\end{eqnarray*}
Since $u_0 \le t_{early}=\frac{1.01}{\sigma(1-\gamma)-\rho} \log (2N)$,
the conclusion in (c) follows.
\qed

\vspace{.2cm}
\noindent
{\em Proof of Proposition~\ref{prop:Z10}(a-b).}
We define $t_2=t_0+t_{early}$ and
\begin{eqnarray*}
	E_{21} &=& \left\{ T_{10;c_{10,3}} \le S_{10,01,rec} \wedge
                        t_2 \right\}, \\
	E_{22} &=& \{ X_{01}(T_{10;c_{10,3}})
		\ge 1-c_{10,3}-(2N)^{-\delta_{01,1}/2} \}, \\
	F_1 &=& \left\{S_{10,01,rec} \le T_{10;c_{10,3}} \wedge
                        (t_2 \vee t'_{01;1-c_{01,2}}) \right\},
\end{eqnarray*}
then $E_{21}\cap E_{22} \subset E_2$. Before
$T_{10;c_{10,3}} \wedge (t_2 \vee t'_{01;1-c_{01,2}})$,
the rate of recombination events between type 10 and 01 individuals
is at most $4\rho N X_{10} X_{01} \le 4\rho N N^{-b_{10,3}}
\le C_{\rho} N^{-b_{10,3}}$. Hence the total number of recombination
events between type 10 and 01 individuals before
$T_{10;c_{10,3}} \wedge (t_2 \vee t'_{01;1-c_{01,2}})$
is dominated by a Poisson random variable
with mean $C_{\rho,\gamma,\sigma} N^{-b_{10,3}} \log N$.
Therefore
\begin{eqnarray}
	\Pbold (F_1) \le C_{\rho,\gamma,\sigma} N^{-7b_{10,3}/8}.
	\label{ineq:F1}
\end{eqnarray}
On $F_1^c$, we have $S_{10,01,rec} > T_{10;c_{10,3}}$ or $S_{10,01,rec} > t_2$,
We observe that 
\begin{eqnarray}
	E_{21}^c \cap F_1^c
	&=&  \{ T_{10;c_{10,3}} > S_{10,01,rec} \vee t_2 \}
		\cap \{ S_{10,01,rec} > T_{10;c_{10,3}} \vee t_2 \} \nnb \\
	&\subset& \{S_{10,01,rec}>t_2, \ T_{10;c_{10,3}} > t_2 \}.
	\label{ineq:2sets}
\end{eqnarray}
Therefore Lemma~\ref{lem:10early}(c) implies
\begin{eqnarray}
	\Pbold (E_{21}^c \cap F_1^c \cap A_4^c \cap B_4)
	\le C_{\rho,\gamma,\sigma} N^{-\delta_{10,1}}. \label{ineq:prp31ab1}
\end{eqnarray} 
Let $F_4=\{T_{10;c_{10,2}} \le t'_{01;1-c_{01,2}} \wedge T_{11;1/(2N)} \}$,
then reasoning similar to that of~(\ref{ineq:2sets}) implies
\begin{eqnarray*}
	F_4^c \cap
		\{ t'_{01;1-c_{01,2}} \ge T_{11;1/(2N)} \vee T_{10;c_{10,2}} \}
	\subset \{T_{10;c_{10,2}} \vee t'_{01;1-c_{01,2}} \ge T_{11;1/(2N)}\},
\end{eqnarray*}
which implies
\begin{eqnarray*}
	\lefteqn{ \Pbold(\{ t'_{01;1-c_{01,2}} \ge T_{11;1/(2N)}
			\vee T_{10;c_{10,2}} \}
		\cap E_{21} \cap F_1^c \cap A_4^c \cap B_4) } \\
	&\le& \Pbold( \{ T_{10;c_{10,2}}
			\vee t'_{01;1-c_{01,2}} \ge T_{11;1/(2N)} \}
		 \cap E_{21} \cap F_1^c \cap A_4^c \cap B_4) \\
	&& \qquad + \Pbold(F_4 \cap A_4^c \cap B_4). \nnb
\end{eqnarray*}
The first set on the right hand side satisfies
\begin{eqnarray*}
	\lefteqn{\{ T_{10;c_{10,2}} \vee t'_{01;1-c_{01,2}} \ge T_{11;1/(2N)}\}
                 \cap E_{21} \cap F_1^c } \\
	&& \subset \{ T_{10;c_{10,2}} \vee t'_{01;1-c_{01,2}} \ge T_{11;1/(2N)} 
		> T_{10;c_{10,3}} \vee t_2 \} \cap E_{21} \\
	&& \subset \{ T_{10;c_{10,2}} \vee t'_{01;1-c_{01,2}} \ge T_{11;1/(2N)}
                > t_2 \ge T_{10;c_{10,3}} \} = \emptyset,
\end{eqnarray*}
therefore
\begin{eqnarray}
        \lefteqn{ \Pbold(\{ t'_{01;1-c_{01,2}} \ge T_{11;1/(2N)}
		\vee T_{10;c_{10,2}} \cap E_{21} \cap F_1^c \cap A_4^c \cap B_4)
		} \nnb \\
        &&\qquad \le \Pbold(F_4 \cap A_4^c \cap B_4)
	\le C_{\rho,\gamma,\sigma} N^{-\delta_{10,1}}
	\label{ineq:prp31ab2}
\end{eqnarray}
by Lemma~\ref{lem:10early}(b).
On $\{ t'_{01;1-c_{01,2}} < T_{11;1/(2N)} \vee T_{10;c_{10,2}} \}$,
we have $T_{10;c_{10,3}}\ge T_{10;c_{10,2}} \ge t'_{01;1-c_{01,2}}$,
therefore Lemma~\ref{lem:01early} implies
\begin{eqnarray}
	\lefteqn{ \Pbold (\{ X_{10}(T_{10;c_{10,3}})+X_{01}(T_{10;c_{10,3}})
		\le 1-(2N)^{-\delta_{01,1}/2} \} } \label{ineq:prp31ab3} \\
	&& \cap \{ t'_{01;1-c_{01,2}} < T_{11;1/(2N)} \vee T_{10;c_{10,2}} \}
	\cap E_{21} \cap F_1^c \cap A_4^c \cap B_4)
	\le C_{\rho,\gamma\sigma} N^{-\delta_{01,1}}. \nnb
\end{eqnarray}
Combining~(\ref{ineq:F1}),~(\ref{ineq:prp31ab1}),~(\ref{ineq:prp31ab2}),
and~(\ref{ineq:prp31ab3}) yields
\begin{eqnarray*}
	\Pbold((E_{22}^c \cup E_{21}^c) \cap A_{41}^c \cap A_{42}^c \cap B_4
		\cap E_1)
	\le C_{\rho,\gamma\sigma} (N^{-\delta_{01,1}}
		+ N^{-\delta_{10,1}} + N^{-7b_{10,3}/8}),
\end{eqnarray*}
where we also recall from Lemma~\ref{lem:10early} that
$A_4^c=A_{41}^c \cap A_{42}^c \cap E_1$.
We further combine the above estimate with the first two statements of
Lemma~\ref{lem:10early}(a) to obtain
\begin{eqnarray}
        \lefteqn{ \Pbold((E_{22}^c \cup E_{21}^c) \cap B_4 \cap E_1) } \nnb \\
        && \le C_{\rho,\gamma\sigma} (N^{-\delta_{01,1}}
                + N^{-\delta_{10,1}} + N^{-7b_{10,3}/8}
		+ \delta_{10,0} t_0 + N^{-a_1}).
	\label{ineq:10mid3}
\end{eqnarray}

It remains to show that $B_4^c =  \{t_0 > T_{01;c_{01,0}} \vee T_{11;1/(2N)}\}$
has a small probability.
Let $F_2 = \{ T_{01;c_{01,0}} < t_0 \wedge T_{11;1/(2N)} \}$.
Before $T_{11;1/(2N)}$, the jump rates of $X_{01}$ satisfy
\[ r^+_{01} \le N (1+\sigma\gamma+\rho) X_{01} (1-X_{01}), \
        r^-_{01} \ge N (1-\sigma\gamma+\rho) X_{01} (1-X_{01}). \]
We take $\check\xi=X_{01}$, $\alpha=1+\rho$, $\theta=\sigma\gamma$,
$\delta_0=0$, $\delta_1=0.9 c_{01,0}$, $\delta_2 = (1-\zeta)/4$, and
$Y(t) = \LL(t;(2N)^{-\zeta},\sigma\gamma)$ in Lemma~\ref{lem:1dim_early}
to obtain
\[ \Pbold\left( X_{01}(s) \ge c_{01,0} \mbox{ for some }
                s\le t_{01;0.9 c_{01,0}} \wedge T_{11;1/(2N)} \right)
        \le C_{\rho,\gamma,\sigma} N^{-(1-\zeta)/4}. \]
By the choice of $a_0$ in~(\ref{eq:constants}),
$t_0 = \frac{\zeta}{3\sigma\gamma}\log (2N)
= \frac{1}{\sigma\gamma}\log (2N)^{\zeta/3}
< t_{01;0.9 (2N)^{-\zeta/3}} = t_{01;0.9 c_{01,0}}$, therefore
\begin{eqnarray*}
        \Pbold(F_2)
        \le \Pbold(T_{01;c_{01,0}} \le t_{01;0.9 c_{01,0}} \wedge T_{11;1/(2N)})
		\le C_{\rho,\gamma,\sigma} N^{-(1-\zeta)/4}.
\end{eqnarray*}
We observe that $B_4^c \cap F_2^c \subset \{t_0 \wedge T_{01;c_{01,1}}
> T_{11;1/(2N)} \}$. By an argument similar to the one leading
to~(\ref{ineq:F1}),
$\Pbold(B_4^c \cap F_2^c) \le C_{\rho,\gamma,\sigma} N^{-\zeta/4}$,
which implies
\begin{eqnarray}
	\Pbold(B_4^c) \le C_{\rho,\gamma,\sigma}
		(N^{-(1-\zeta)/4} + N^{-\zeta/4}). \label{ineq:B4}
\end{eqnarray}
Combining~(\ref{ineq:10mid3}) and~(\ref{ineq:B4}) yields the desired result in
(b). For part (a), we combine the third statement of
Lemma~\ref{lem:10early}(a) and~(\ref{ineq:B4}) to obtain
the desired result.
\qed

\vspace{.2cm}
\noindent
{\em Proof of Proposition~\ref{prop:Z10}(c-e).}
Recall that $Z_{10}(T_{10;c_{10,3}}+t)=\LL(t; c_{10,3},\sigma(1-\gamma))$
for $t\in [T_{10;c_{10,3}},T_{Z_{10};1-c_{10,3}}]$, and
$T_{Z_{10};1-c_{10,3}} = T_{10;c_{10,3}} +\frac 1 {
                \sigma(1-\gamma)}\log \frac{1-c_{10,3}}{c_{10,3}}$.
We work on $t\ge T_{10;c_{10,3}}$ throughout this proof.
On $E_2\cap E_1$, we have
$X_{01}(T_{10;c_{10,3}})\ge 1-c_{10,3}-(2N)^{-\delta_{01,1}/2}$,
$X_{10}(T_{10;c_{10,3}})=c_{10,3}$ and $X_{00}(T_{10;c_{10,3}})
\le (2N)^{-\delta_{01,1}/2}$.
We can then write down the following equation
using the jump rates of $X_{10}$ in~(\ref{eq:rij}):
\begin{eqnarray*}
        X_{10}(t) &=& c_{10,3} + M_{10}(t) + \int_{T_{10;c_{10,3}}}^t
                X_{10}(s) [\sigma(1-\gamma) (1-X_{10}(s)) \\
	&& - \sigma X_{11}(s) + \sigma\gamma X_{00}(s)]
        	+ \rho (X_{11}(s) X_{00}(s)-X_{10}(s) X_{01}(s)) \ ds,
\end{eqnarray*}
where $M_{10}$ is a martingale with maximum jump size $1/(2N)$ and
quadratic variation
$ \langle M_{10} \rangle(t) = \frac{1}{2N} \int_{T_{10;c_{10,3}}}^t
	(1+\rho) X_{10}(s)(1-X_{10}(s)) + \rho X_{11}(s) X_{00}(s) \ ds$.
We use Lemma~\ref{lem:1dim_mid} 
with $\theta=\sigma(1-\gamma)$, $u_1=0$,
$u_2=\frac 1 {\sigma(1-\gamma)}\log \frac{1-c_{10,3}}{c_{10,3}}$,
$\delta_1=\delta_{01,1}/4$, $\delta_2=\infty$,
$\epsilon_0=\epsilon_1=b_{10,3}=\delta_{01,1}/10$,
$T=T_{11;\delta_{11}} \wedge T_{00;(2N)^{-\delta_1}}$,
$\epsilon_2(t)= - \sigma X_{11}(t) + \sigma\gamma X_{00}(t)$,
$\epsilon_3(t)=\rho (X_{11}(t)X_{00}(t) - X_{01}(t)X_{10}(t))$,
$\epsilon_4(t)=X_{11}(t)X_{00}(t)$,
$Y(t)=Z_{10}(T_{10;c_{10,3}}+t)$,
and $D_1=E_2\cap E_1$ to obtain
\begin{eqnarray}
        \lefteqn{ \Pbold \left( |X_{10}(s,\omega)-Z_{10}(s,\omega)|
                > (2N)^{-\delta_{10,2}} \mbox{ for some } \omega\in E_2\cap E_1,
		\right. } \label{ineq:10mid1} \\
        && \left. s \in [T_{10;c_{10,3}}, T_{Z_{10};1-c_{10,3}} \wedge
		T_{11;\delta_{11}} \wedge T_{00;(2N)^{-\delta_{01,1}/4}}]
		\right) \le (2N)^{-\delta_{10,2}}, \nnb
\end{eqnarray}
where $\delta_{10,2}=(\delta_1-\epsilon_1-\epsilon_2)/3=\delta_{01,1}/60$,
as defined in~(\ref{eq:constants}).
The jump rates of $X_{00}$ satisfy
\[ r^+_{00} \le N [(1-\sigma\gamma+\rho) X_{00} (1-X_{00})
	+ 2 \rho X_{01}X_{10}],
	r^-_{00} \ge N (1+\sigma\gamma+\rho) X_{00} (1-X_{00}). \]
On $E_2\cap E_1$, we have $X_{00}(T_{10;c_{10,3}})\le (2N)^{-\delta_{01,1}/2}$.
Therefore by Lemma~\ref{lem:1dim_late},
\begin{eqnarray*}
	\Pbold\left( \left\{
		\sup_{s\in [T_{10;c_{10,3}},T_{Z_{10};1-c_{10,3}}]} X_{00}(s)
		\ge (2N)^{-\delta_{01,1}/4} \right\} \cap E_2\cap E_1 \right)
	\le C N^{-1/2}.
\end{eqnarray*}
We combine the above and~(\ref{ineq:10mid1}) to arrive at the
desired conclusion of (c).

For (d), we observe that
the jump rates of $X_{11+10}=X_{11}+X_{10}$ satisfy
\begin{eqnarray*}
        r^+_{11+10} &=&  N X_{11} [(1+\sigma+\rho) X_{01}
                        + (1+\sigma(1+\gamma)+2\rho) X_{00}] \\
        &&      + N X_{10} [(1+\sigma(1-\gamma)+2\rho) X_{01}
                        + (1+\sigma+\rho) X_{00}] \\
        r^-_{11+10} &=& N X_{11} [(1-\sigma+\rho) X_{01}
                        + (1-\sigma(1+\gamma)+\rho) X_{00}] \\
        && + N X_{10} [(1-\sigma(1-\gamma)+\rho) X_{01}
                        + (1-\sigma+\rho) X_{00}],
\end{eqnarray*} 
where we drop the terms involving $X_{11} X_{10}$ in $r^\pm_{10}$ and
$r^\pm_{11}$, which correspond to type 11 individuals replaced by
type 10 individuals or {\em vice versa}.
Therefore $X_{11+10}$ dominates $1-\eta$ where we define $\eta$ to be a
jump process with initial condition
$\eta(T_{10;c_{10,3}})=1-X_{11+10}(T_{10;c_{10,3}})$ and jump rates of
\begin{eqnarray*}
        r^+_{\eta} = N (1-\sigma(1-\gamma)+\rho) \eta (1-\eta), \
        r^-_{\eta} = N (1+\sigma(1-\gamma)+\rho) \eta (1-\eta).
\end{eqnarray*}
Since $\eta(T_{Z_{10};1-c_{10,3}})\le 1-X_{10}(T_{Z_{10};1-c_{10,3}})
\le c_{10,3}$ on $E_4 \cap E_3\cap E_2\cap E_1$,
by Lemma~\ref{lem:1dim_late},
\begin{eqnarray*}
        \Pbold \left( \{ \eta(t) \ge \sqrt{c_{10,3}}
                \mbox{ for some } t \ge T_{Z_{10};1-c_{10,3}} \}
        \cap E_4 \cap E_3\cap E_2\cap E_1 \right) \le C N^{-1/2}.
\end{eqnarray*}
This implies the desired conclusion of (d).

Let $\tilde\eta$ be a time change of $\eta$ by $1-\eta$, then
$2N \tilde\eta$ is a branching process and the clock for $\tilde\eta$ runs
at the rate of at most 1.02 times that of $\eta$
on $\{\tilde\eta(t) < \sqrt{c_{10,3}}
\mbox{ for all } t \ge T_{Z_{10};1-c_{10,3}} \} \cap E_4 \cap E_3\cap E_2\cap E_1$.
By Lemma~\ref{lem:branching}(b),
\begin{eqnarray*}
	P(\{ \tilde\eta\left(T_{Z_{10};1-c_{10,3}}+0.99t_{late} ) >0 \}
	\cap E_4 \cap E_3\cap E_2\cap E_1 \right)
	\le C N c_{10,3} e^{-\log (2N)}.
\end{eqnarray*}
Hence
$P\left(\left\{ \eta\left(T_{Z_{10};1-c_{10,3}}
                +t_{late} \right) >0 \right\}
        \cap E_4 \cap E_3\cap E_2\cap E_1 \right) \le C c_{10,3}$,
which implies (e) since
$T_{Z_{10};1-c_{10,3}} + t_{late} = T_\infty$.
\qed

\vspace{.2cm}
\noindent
{\em Proof of Proposition~\ref{prop:Z10}(g-h).}
We define $c_4$ and $\delta_{10,4}$ such that
$c_4=\max(\sqrt{c_1}+\delta_{11},c_2,c_3)\le N^{-2\delta_{10,4}}$ and we let
\[ S_{X,Z,far} = \inf\{ t\ge T_{10;c_1}:
        |X_{10}(t)-Z_{10}(t)| \vee X_{00}(t)> c_4 \}. \]
By Proposition~\ref{prop:Z10}(c,d), there exists $\delta_{10,3}>0$ such that
\begin{eqnarray}
        \lefteqn{ \Pbold( \{ S_{X,Z,far}\le T_{11;\delta_{11}} \}
		\cap E_2\cap E_1) } \nnb \\
	&& \le \Pbold( (E_3^c \cup (E_5^c\cap E_4\cap E_3)) \cap E_2\cap E_1)
        \le C_{\rho,\gamma,\sigma} N^{-\delta_{10,3}}, \label{ineq:SXZ11}
\end{eqnarray}
where we have used that on $E_4\cap E_3$, $S_{X,Z,far} \ge T_{Z_{10};1-c_1}$
and on $E_5$,
$X_{10}(t) > 1-\sqrt{c_1}-X_{11}(t) > 1-\sqrt{c_1}-\delta_{11}$
and $X_{00}(t) \le 1-X_{10}(t)-X_{11}(t) < \sqrt{c_1}$ for
$t\ge T_{Z_{10};1-c_1}$. Notice that on $E_2\cap E_1$,
$X_{11}(t)=0=Z_{11}(t)$ for all $t\le T_{10;c_1}$.
For $t<S_{X_{11},Z_{11},diff} \wedge S_{X,Z,far} \wedge T_{11;\delta_{11}}$,
we have
\begin{eqnarray*}
        |r^+_{Z,11}-r^+_{11}|
        \le N \delta_{11} [(\sigma-\rho) 3c_4+\delta_{11}]
                + 2\rho N (3c_4+\delta_{11})
        \le 4 N \delta_{11} c_4
\end{eqnarray*}
and similarly, $|r^-_{Z,11}-r^-_{11}| \le 4 N \delta_{11} c_4$.
Thus the absolute difference between $X_{11}$ and $Z_{11}$ is bounded
above by a Poisson process of rate $8 N \delta_{11} c_4$, which stays 0 during
$[T_{10;c_1},T_{10;c_1}+t_{mid}+t_{late}]$ with probability at least 
$1-c_4^{1/2}$, if $t_{mid}+t_{late}\le c_4^{-1/4}$, which is satisfied by
our choice of $t_{mid}+t_{late}=\Ocal(\log N)$. Hence
\[ \Pbold( \{ S_{X_{11},Z_{11},diff} \le T_\infty \wedge
        S_{X,Z,far} \wedge T_{11;\delta_{11}} \} \cap E_2\cap E_1)
        \le c_4^{1/2}.
\]
We combine~(\ref{ineq:SXZ11}) and the above estimate to obtain
\begin{eqnarray*}
        \Pbold( \{ S_{X_{11},Z_{11},diff} \le T_\infty \wedge
                T_{11;\delta_{11}} \} \cap E_2\cap E_1) 
        \le c_4^{1/2} + C_{\rho,\gamma,\sigma} N^{-\delta_{10,3}},
\end{eqnarray*}
which implies (g).

Let $F_3=\{ T_{Z_{11};\{0,\delta_{11} \}} \ge T_{Z_{10};1-c_1} \}$.
Starting from $T_{Z_{10};1-c_1}$, $Z_{11}$ is a time-changed branching 
process. We perform a time change of $1-Z_{11}$ (from time $T_{Z_{10};1-c_1}$
onwards) to obtain a branching process $\tilde Z_{11}$,
then the clock for $\tilde Z_{11}$ runs faster than that of $Z_{11}$
(at a rate of at most $1/(1-\delta_{11})$ times
before $\tilde Z_{11}$ reaches $\delta_{11}$).
From time $T_{Z_{10};1-c_1}$ onwards, 0 and $\delta_{11}$
are absorption points for
$Z_{11}(\cdot \wedge T_{Z_{11};\delta_{11}})$,
We use Lemma~\ref{lem:branching}(d) below to deduce that
\begin{eqnarray*}
        \lefteqn{ \Pbold\left(\{ Z_{11}(T_\infty \wedge
                        T_{Z_{11};\delta_{11}}) \in (0,\delta_{11}) \}
                \cap F_3 \cap E_2\cap E_1 \right) } \nnb \\
        &\le& \Pbold\left(\{ \tilde Z_{11}(s) \in (0,\delta_{11})
                \mbox{ for all } s\le (1-\delta_{11})T_\infty \}
                \cap F_3 \cap E_2\cap E_1 \right) \nnb \\
        &\le& (2N\delta_{11})^2 C_{\rho,\gamma,\sigma} \exp(-0.99
		\sigma\gamma (T_\infty-T_{Z_{10};1-c_1})), \nnb \\
	&\le& C_{\rho,\gamma,\sigma} (\log^2 N) \exp(-0.99 \sigma\gamma t_{late})
	\le C_{\rho,\gamma,\sigma} N^{-\delta_{10,4}},
\end{eqnarray*}
if we choose a sufficiently small $\delta_{10,4}$. Therefore
\begin{eqnarray*}
        \Pbold\left(\{ Z_{11}(T_\infty) \in (0,\delta_{11}), \
                        T_{Z_{11};\delta_{11}}> T_\infty \}
                \cap F_3 \cap E_2\cap E_1 \right)
        \le C_{\rho,\gamma,\sigma} N^{-\delta_{10,4}}.
\end{eqnarray*}
On $\{ S_{X_{11},Z_{11},diff} > T_\infty \wedge T_{11;\delta_{11}} \}$,
$X_{11}$ and $Z_{11}$ agree up to $T_\infty \wedge T_{11;\delta_{11}}$.
Therefore
\begin{eqnarray*}
        \lefteqn{ \Pbold \left( \{ S_{X_{11},Z_{11},diff} > T_\infty,
                        T_{11;\delta_{11}} > T_\infty,
                        X_{11}(T_\infty)=Z_{11}(T_\infty)
                                \in (0,\delta_{11}), \right. } \\
        && \qquad\qquad\qquad\qquad \left.
                        T_{11;\{0,\delta_{11} \}}
                                \ge T_{Z_{10};1-c_1} \}
                \cap E_2\cap E_1 \right)
        \le C_{\rho,\gamma,\sigma} N^{-\delta_{10,4}}.
\end{eqnarray*}
We can drop the condition
$T_{11;\{0,\delta_{11} \}} \ge T_{Z_{10};1-c_1}$, since on
$\{ S_{X_{11},Z_{11},diff} > T_\infty, T_{11;\delta_{11}} > T_\infty,
        T_{11;\{0,\delta_{11} \}} < T_{Z_{10};1-c_1} \}$,
we have $X_{11}(T_\infty)=Z_{11}(T_\infty)=0$.
Hence
\begin{eqnarray*}
        \Pbold \left( \{ T_{11;\delta_{11}} > T_\infty,
                        X_{11}(T_\infty)=Z_{11}(T_\infty)
                                \in (0,\delta_{11}) \}
                \cap E_7 \cap E_2\cap E_1 \right)
        \le C_{\rho,\gamma,\sigma} N^{-\delta_{10,4}},
\end{eqnarray*}
which implies the desired result in (h).
\qed

\section{Supporting Lemmas}
\label{sec:lemmas}
In this section, we establish
Lemmas~\ref{lem:1dim_early} to~\ref{lem:1dim_late}, one each
for the early, middle, and late phase. They are used for the proof of
Proposition~\ref{prop:Z10} in \S\ref{sec:case1}.
Lemma~\ref{lem:1dim_early} deals with the early
phase and approximates a 1-dimensional jump process undergoing selection by
a deterministic function, where the error bound depends only on the
initial condition of the process, as long as the process is stopped
before it reaches $\Ocal(1)$. Lemma~\ref{lem:1dim_mid} deals with the middle
phase and uses the logistic growth as an approximation. The main
difference between the early phase and the middle phase is the error bound:
in Lemma~\ref{lem:1dim_mid}, the error bound depends on both the initial
and terminal conditions of the process. Lemma~\ref{lem:1dim_late} deals
with the late phase, for which we only need to show that the process
does not stray too far away from 1 (or 0 for $X_{00}$) once it gets close to
1 (or 0).

\begin{lemma}
Let $\alpha\ge 1$, $\theta\in (0,1)$, $\delta_0\in [0,1/2]$
and $x\in (0,1]$ be constants.
Let $\xi$ be a jump process with initial value $\xi(0)=(2N)^{-x}\ge (2N)^{-1}$,
jump size $1/(2N)$, and jump rates
\begin{eqnarray*}
	r^+ = N \xi [(\alpha+\theta) (1-\xi)-\delta_0], \
	r^- = N \xi [(\alpha-\theta) (1-\xi)+\delta_0].
\end{eqnarray*}
Suppose $Y$ is a deterministic process that satisfies
\[ Y(t) = (2N)^{-x} + \int_0^t Y(s) (\theta(1-Y(s))-\delta_0) \ ds. \]
If $u_0=\inf\{t: Y(t)=\delta_1\}
\le (\log 2)/(3\theta\delta_1+\delta_0)$, then
there exists $\delta_2 \in (0,(1-x)/4]$ such that
\[ P \left(|\xi(s)-Y(s)| > 4 N^{-\delta_2} Y(s)
		\mbox{ for some } s\le u_0 \right)
	\le C_{\alpha,\theta} N^{-\delta_2}.
\]
Moreover, if $\check\xi$ and $\hat\xi$ are jump processes such that
$\hat\xi \ge \xi \ge \check\xi$ before a stopping time $T$, then
$P \left(\hat\xi(s) < (1-4 N^{-\delta_2}) Y(s) \mbox{ for some }
s\le u_0 \wedge T\right) \le C_{\alpha,\theta} N^{-\delta_2}$
and
$P \left(\check\xi(s) > (1+4 N^{-\delta_2}) Y(s) \mbox{ for some }
s\le u_0 \wedge T\right) \le C_{\alpha,\theta} N^{-\delta_2}$.
\label{lem:1dim_early}
\end{lemma}
\proof
We can write
\begin{eqnarray*}
        d\xi = dM_{\xi} + \xi (\theta(1-\xi)-\delta_0) \ dt, \
        d \langle M_{\xi} \rangle = \frac \alpha {2N} \xi (1-\xi) \ dt,
\end{eqnarray*}
and consequently,
\begin{eqnarray}
        d (e^{-\theta t} \xi(t)) &=& d\tilde M_{\xi}(t)
		- e^{-\theta t} (\theta\xi(t)^2+\delta_0\xi(t)) \ dt
	\label{eq:tildexi} \\
        d\langle \tilde M_{\xi} \rangle(t) &=& \frac \alpha {2N}
                e^{-2\theta t} \xi(t) (1-\xi(t)) \ dt. \nnb
\end{eqnarray}
We define $\tau = \inf\{t\le u_0: \xi(t)\ge 2\delta_1 \}$,
and take expectation on both sides of~(\ref{eq:tildexi}) to obtain
\begin{eqnarray*}
	E [ e^{-\theta (t\wedge\tau)} \xi(t\wedge\tau) ] = (2N)^{-x}
		- E\left[\int_0^{t\wedge\tau} e^{-\theta s}
			(\theta\xi(s)^2 + \delta_0\xi(s))
			\ ds \right]
	\le (2N)^{-x}.
\end{eqnarray*}
As in the steps leading to~(\ref{ineq:JB}), we use
Jensen's and Burkholder's inequalities to obtain
\begin{eqnarray}
        \lefteqn{ E \left[ \sup_{s\le t\wedge\tau} |\tilde M_{\xi}(s)| \right]
	\le \frac C N + \frac{C_\alpha}{N^{1/2}} \left( E \left[ \int_0^t
                e^{-2\theta s} \xi(s) {\bf 1}_{\{ s\le\tau \} } \ ds \right]
		\right)^{1/2} } \nnb \\
	&&\le \frac C N + \frac{C_\alpha}{N^{1/2}} \left( \int_0^t
                e^{-\theta s} (2N)^{-x} \ ds \right)^{1/2}
        \le C_{\alpha,\theta} N^{-(1+x)/2}.
	\label{ineq:tildeMxi}
\end{eqnarray}
Since $de^{-\theta t} Y(t)
= - e^{-\theta t} (\theta Y(t)^2+\delta_0 Y(t)) \ dt$,
we use~(\ref{ineq:tildeMxi}) in~(\ref{eq:tildexi}) to obtain
\begin{eqnarray*}
	\lefteqn{ E \left[ \sup_{s\le t\wedge\tau} e^{-\theta s}
		|\xi(s)-Y(s)| \right] } \\
	&\le& C_{\alpha,\theta} N^{-(1+x)/2}
        	+ E \left[ \int_0^{t\wedge\tau} e^{-\theta s}
			(\theta |\xi(s)^2-Y(s)^2|+\delta_0 |\xi(s)-Y(s)|
		\ ds \right] \\
	&\le& C_{\alpha,\theta} N^{-(1+x)/2}
                + E \left[ \int_0^t (3\theta\delta_1+\delta_0)
		e^{-\theta s} |\xi(s)-Y(s)| {\bf 1}_{ \{s\le\tau\} } \ ds \right] \\
	&\le& C_{\alpha,\theta} N^{-(1+x)/2}
                + \int_0^t (3\theta\delta_1+\delta_0)
		E\left[ \sup_{s\le s'\wedge\tau} e^{-\theta s} |\xi(s)-Y(s)|
		\right] \ ds'.
\end{eqnarray*}
Gronwall's inequality implies
\begin{eqnarray*}
	E \left[ \sup_{s\le t\wedge\tau} e^{-\theta s} |\xi(s)-Y(s)|
	\right] \le C_{\alpha,\theta} N^{-(1+x)/2}
		e^{(3\theta\delta_1+\delta_0) t}
	\le C_{\alpha,\theta} N^{-(1+x)/2},
\end{eqnarray*}
since $\tau\le u_0\le (\log 2)/(3\theta\delta_1+\delta_0)$.
Let $\delta_2\in (0,(1-x)/4]$, then
\[ P \left(|\xi(s)-Y(s)| \ge N^{-\delta_2-x} e^{\theta s} \mbox{ for some }
	s \le u_0\wedge\tau \right)  \le C_{\alpha,\theta} N^{-\delta_2}. \]
We observe that for $s\le u_0$,
$(2N)^{-x} e^{(\theta-\theta\delta_1-\delta_0) s} \le Y(s)$, hence
$N^{-x} e^{\theta s}/Y(s) \le 2^x e^{(\theta\delta_1+\delta_0)s}
\le 2^x e^{(\theta\delta_1+\delta_0)(\log 2)/(3\theta\delta_1+\delta_0)}
\le 4$, i.e. $N^{-x} e^{\theta s} \le 4 Y(s)$. Hence
\[ P \left(|\xi(s)-Y(s)| \ge 4 N^{-\delta_2} Y(s) \mbox{ for some }
        s \le u_0\wedge\tau \right)  \le C_{\alpha,\theta} N^{-\delta_2}. \]
We can drop $\tau$ in the event above, since $|\xi(\tau)-Y(\tau)|\ge Y(\tau)$.
The conclusion follows.
\qed

\begin{lemma}
Let $\theta,\epsilon_0,\epsilon_1 \in (0,1)$ and
$a_0,a_1>0$ be constants. Suppose
$Y$ is a deterministic process defined from a stopping time $u_1$ onwards
that has initial condition $Y(u_1)=b_0\ge a_0 (2N)^{-\epsilon_0}$ and satisfies
\[ Y(t) = b_0 + \int_{u_1}^t \theta Y(s) (1-Y(s)) \ ds. \]
Let $u_2=u_1 +\frac 1 {\theta} \log \frac{1-b_0}{b_0} \frac{1-b_1}{b_1}$
such that $Y(u_2)=1-b_1 \le 1-a_1 (2N)^{-\epsilon_1}$.
Suppose $T$ is a stopping time and
$\xi$ is jump process that takes values in $[0,1]$,
has jump size $1/(2N)$ and satisfies
\begin{eqnarray*}
        \xi(t\wedge T) &=& \xi(u_1) + M(t\wedge T) + \int_{u_1}^{t\wedge T}
		\xi(s) [\theta(1-\xi(s)) + \epsilon_2(s)]
			 + \epsilon_3(s) \ ds \\
        \langle M \rangle (t\wedge T) &=& \frac{1+\rho}{2N}
                \int_{u_1}^{t\wedge T} \xi(s)(1-\xi(s)) + \epsilon_4(s) \ ds,
\end{eqnarray*}
where $|\epsilon_2(t)|,|\epsilon_3(t)|\le (2N)^{-\delta_1}$,
$\epsilon_4(t)\le 1$ for $t\le T$,
and $M$ is a jump martingale with jump size $1/(2N)$.
Furthermore, suppose on a set $D_1\in \Fcal(u_1)$, we have
$|\xi(u_1)-b_0| \le (2N)^{-\delta_2}$.
We define $D_2 = \{ T \ge u_1 \}$ and $\delta_3$ to be a constant
$\le ((\delta_1\wedge \delta_2\wedge \frac 1 2)-\epsilon_0-\epsilon_1)/3$.
If $\delta_3>0$, then
\begin{eqnarray*}
        P \left( \left\{
		\sup_{s \in [u_1,u_2\wedge T]} |\xi(s,\omega)-Y(s,\omega)|
                > (2N)^{-\delta_3} \right\} \cap D_1\cap D_2\right)
        \le (2N)^{-\delta_3}.
\end{eqnarray*}
\label{lem:1dim_mid}
\end{lemma}
\proof
Let $D=D_1\cap D_2$. Notice that $D\in \Fcal(u_1)$. Since
\begin{eqnarray*}
        \lefteqn{ |\xi(t) [\theta(1-\xi(t)) + \epsilon_2(t)]
                - \theta Y(t) (1-Y(t))| {\bf 1}_{\{t\le T\} } } \\
        && \le (2N)^{-\delta_1}+\theta |\xi(t)-Y(t)||1-\xi(t)-Y(t)| {\bf 1}_{\{t\le T\} }
		\\
        && \le (2N)^{-\delta_1} + \theta |\xi(t)-Y(t)| {\bf 1}_{\{t\le T\} },
\end{eqnarray*}
we have
\begin{eqnarray*}
        \lefteqn{ |\xi((u_1+t)\wedge T)-Y((u_1+t)\wedge T)| {\bf 1}_D
		\le |\xi(u_1) - Y(u_1)| } \\
	&& \qquad + |M((u_1+t)\wedge T) {\bf 1}_D| + \int_{u_1}^{(u_1+t)\wedge T}
		 	[(2N)^{-\delta_1} + \theta |\xi(s)-Y(s)|] {\bf 1}_D \ ds.
\end{eqnarray*}
By Jensen's and Burkholder's inequalities, \nocite{Burkholder73}
\begin{eqnarray*}
        E\left[ \sup_{u_1\le s\le u_1+t} |M(s\wedge T) {\bf 1}_D| \right]
        \le \frac C N + C \sqrt {\frac t N}
        \le C \sqrt{\frac t N},
\end{eqnarray*}
therefore
\begin{eqnarray*}
        \lefteqn{ E\left[ \sup_{u_1\le s\le u_1+t}
                |\xi(s\wedge T)-Y(s\wedge T)| {\bf 1}_D \right]
        \le C \sqrt{\frac t N} + E \left[ |\xi(u_1) - Y(u_1)| {\bf 1}_D \right]
	} \\
        && \qquad \qquad \qquad \qquad + 2 (2N)^{-\delta_1} t
        	+ \int_{u_1}^{u_1+t} \theta E[|\xi(s)-Y(s)| {\bf 1}_{\{s\le T\} } {\bf 1}_D] \ ds.
\end{eqnarray*}
Since $E[|\xi(s)-Y(s)| {\bf 1}_{\{s\le T\} } {\bf 1}_D]
\le E[|\xi(s\wedge T)-Y(s\wedge T)| {\bf 1}_D]$, and
$|\xi(u_1) - Y(u_1)| {\bf 1}_D \le (2N)^{-\delta_2}$, we have
\begin{eqnarray*}
        E\left[ \sup_{u_1\le s\le u_1+t}
                |\xi(s\wedge T)-Y(s\wedge T)| {\bf 1}_D \right]
        \le C \left(\sqrt{\frac t N} + \frac 1 {(2N)^{\delta_2}}
		+ \frac{t}{(2N)^{\delta_1}} \right) e^{\theta t}
\end{eqnarray*}
by Gronwall's inequality. We observe that
$u_2-u_1 \le \frac 1 {\theta} \log \frac{1}{b_0 b_1} \le \frac 1 {\theta}
	[(\epsilon_0+\epsilon_1) \log (2N) - \log (a_0 a_1)]$,
therefore the estimate above implies
\begin{eqnarray*}
        E\left[ \sup_{u_1\le s\le u_2}
                |\xi(s\wedge T)-Y(s\wedge T)| {\bf 1}_D \right]
	\le \frac{C (2N)^{\epsilon_0+\epsilon_1} \log N}{a_0 a_1}
		(2N)^{-(\delta_1 \wedge \delta_2 \wedge \frac 1 2)}.
\end{eqnarray*}
Since $0<\delta_3 \le ((\delta_1\wedge \delta_2\wedge \frac 1 2)
	-\epsilon_0-\epsilon_1)/3$, we have
\begin{eqnarray*}
        E\left[ \sup_{u_1\le s\le u_2}
                |\xi(s\wedge T)-Y(s\wedge T)| {\bf 1}_D \right]
        \le (2N)^{-2\delta_3},
\end{eqnarray*}
which implies the desired conclusion.
\qed

\begin{lemma}
Let $\alpha\ge 1$, $\theta\in (0,1)$,
$x\in (0,1]$, $c>0$ and $\kappa\ge 0$ be constants.
Let $\eta\le\hat\eta$ be jump processes where $\eta$ has
initial value $\eta(0)=1-c (2N)^{-x}$, jump size $1/(2N)$, jump rates
\begin{eqnarray*}
        r^+ = N (\alpha+\theta) \eta (1-\eta), \
        r^- = N (\alpha-\theta) \eta (1-\eta) + N\kappa,
\end{eqnarray*}
and absorbing boundary at $1/2$.
For $t\le c (2N)^{-x}/\kappa$ (if $\kappa=0$, then $t=\infty$), we have
\[ P\left( \inf_{s\le t} \hat\eta(s) > 1-(2N)^{-x/2} \right)
\ge P\left( \inf_{s\le t} \eta(s) > 1-(2N)^{-x/2} \right) \ge 1-C N^{-1/2}.
\]
\label{lem:1dim_late}
\end{lemma}
\proof
We take $\xi=1-\eta$ and perform a time change of $1-\xi$ on $\xi$
to obtain a process $\tilde\xi$ with jump rates
\begin{eqnarray*}
        \tilde r^+ = N (\alpha-\theta) \tilde\xi + N\kappa/(1-\tilde\xi), \
        \tilde r^- = N (\alpha+\theta) \tilde\xi.
\end{eqnarray*}
Let $\tilde\xi_{up}$ be a jump process with initial condition
$\tilde\xi_{up}(0)=\tilde\xi(0)=c(2N)^{-x}$, jump size $1/(2N)$ and
jump rates
\begin{eqnarray*}
        \tilde r^+_{up} = N (\alpha-\theta) \tilde\xi_{up} + 2N\kappa, \
        \tilde r^-_{up} = N (\alpha+\theta) \tilde\xi_{up}.
\end{eqnarray*}
Before the stopping time $\tau = \inf\{t\ge 0: \tilde\xi_{up}\ge 1/2 \}$,
$\tilde\xi_{up}$ dominates $\tilde\xi$. We can write
\begin{eqnarray*}
        d\tilde\xi_{up}(t) = dM_{\tilde\xi_{up}}
		+ (\kappa - \theta \tilde\xi_{up}) \ dt,
        d\langle M_{\tilde\xi_{up}} \rangle = \frac 1 {2N}
                (\kappa + \alpha \tilde\xi_{up}(t)) \ dt.
\end{eqnarray*}
Hence $E[\tilde\xi_{up}(t)]=\frac{\kappa}{\theta}
+ \left(c (2N)^{-x}-\frac{\kappa}{\theta} \right) e^{-\theta t}$ and
by Jensen's and Burkholder's inequalities,
\begin{eqnarray*}
	E \left[ \sup_{s\le 2t} M_{\tilde\xi_{up}}(s) \right]
	&\le& \frac C N + \frac C {\sqrt N} \left( \kappa t
		+ \alpha \int_0^{2t} E[\tilde\xi_{up}(t)] \ ds \right)^{1/2}
		\\
	&\le& \frac{C_{\alpha,\theta}}{\sqrt N} \left(
		\kappa t + c (2N)^{-x} \right)^{1/2}
	\le C_{\alpha,\theta} N^{-(1+x)/2},
\end{eqnarray*}
if $\kappa t\le c (2N)^{-x}$, in which case
\begin{eqnarray*}
	P\left( \sup_{s\le 2t} \tilde\xi_{up}(s) \ge (2N)^{-x/2} \right)
	&\le& P\left( \sup_{s\le 2t} M_{\tilde\xi}(s)
                \ge (2N)^{-x/2}-c (2N)^{-x}-4\kappa t \right) \\
	&\le& \frac{C_{\alpha,\theta} N^{-(x+1)/2}}{(2N)^{-x/2}-c(2N)^{-x}-4\kappa t}
	\le C_{\alpha,\theta} N^{-1/2}.
\end{eqnarray*}
On the set $\{\sup_{s\le 2t} \tilde\xi_{up}(s) \le (2N)^{-x/2} \}$,
$\tilde\xi_{up}$ certainly does not reach $1/2$ before time $2t$.
Hence $\tilde\xi_{up}$ dominates $\tilde\xi$ before $2t$
for $\omega\in\{\sup_{s\le 2t} \tilde\xi_{up}(s) \le (2N)^{-x/2} \}$, which
implies $P\left( \sup_{s\le 2t} \tilde\xi(s) < (2N)^{-x/2} \right)
	\ge 1-C_{\alpha,\theta} N^{-1/2}$.
Because $\tilde\xi$ is the process $\xi$ after a time change of $1-\xi$,
the clock for $\tilde\xi$ runs faster than that of $\xi$, but at most twice
as fast before $\tilde\xi$ reaches $1/2$. Therefore
the estimate above implies
        $P\left( \sup_{s\le t} \xi(s) < (2N)^{-x/2} \right)
	\ge P\left( \sup_{s\le 2t} \tilde\xi(s) < (2N)^{-x/2} \right)
	\ge 1-C_{\alpha,\theta} N^{-1/2}$.
The conclusion follows.
\qed

\section{Appendix: A Result on Branching Processes}
\label{sec:aux}
\begin{lemma}
Let $\xi^{(k)}$ be a branching process with $\xi(0)=k$ and
$u(s)=a s^2 + b$ be the probability generating function of the
offspring distribution. Then
\[ G(s,t) = E(s^{\xi^{(k)}(t)}) = \left( \frac{b(s-1)-(as-b)e^{-(a-b)t}}{
        a(s-1)-(as-b)e^{-(a-b)t}} \right)^k. \]
\noindent (a) If $k=1$ and $a>b$, then
\begin{enumerate}
\item $|P(\xi^{(1)}(t)=0)-b/a| \le b e^{-(a-b)t}/a$.
\item $P(1\le \xi^{(1)}(t) \le K) \le C_{a,b} K e^{-(a-b)t}$
	if $K \le e^{(a-b)t}/6$.
\item $P\left(\sup_{s\le t} \xi^{(1)}(s)\ge K\right) \le C_{a,b} e^{(a-b)t}/K$.
\end{enumerate}

\noindent (b) If $a<b$, then $P(\xi^{(k)}(t)>0) \le 1.2 k e^{-(b-a)t}$.

\noindent (c) If $a>b$ and $k\in [1,K]$, then
	$P(\xi^{(k)}(t)\in [1,K])\le k C_{a,b} K e^{-(a-b)t}$.

\noindent (d) If $a>b$ and $\xi$ is a branching process with an initial
	condition that has support on $[0,k]$, then
	$P(\xi(t)\in [1,K])\le k C_{a,b} K e^{-(a-b)t}$. Consequently,
	\[ P(\xi(s)\in [1,K] \mbox{ for all } s\le t)
		\le k C_{a,b} K e^{-(a-b)t}. \]
\label{lem:branching}
\end{lemma}

\proof
The formula for $G(s,t)$ comes from Chapter III.5 of Athreya \& Ney (1972).
\nocite{AthreyaNey}
From this formula, we deduce that
\begin{eqnarray}
	P(\xi^{(k)}(t)=0) = G(0,t)
	= \left( \frac{b-b e^{-(a-b)t}}{a-b e^{-(a-b)t}} \right)^k.
	\label{eq:branch1}
\end{eqnarray}
For (a), we specialise to the case of $k=1$ and $a>b$. We write
$\xi=\xi^{(1)}$, then
\begin{eqnarray*}
	\left| P(\xi(t)=0) - \frac b a \right|
	= \frac{(a-b) be^{-(a-b)t}}{a(a-b e^{-(a-b)t})}
	\le \frac b a e^{-(a-b)t},
\end{eqnarray*}
as required by (a.1). For $s\le 1$, we have
\begin{eqnarray*}
	\lefteqn{ P(1 \le \xi(t) \le K)
	\le s^{-K} \sum_{i=1}^\infty P(\xi(t) = i) s^i
	= s^{-K} (G(s,t)-G(0,t)) } \\
	&&\qquad = \frac{(a-b)^2 s}{(a - b e^{-(a-b)t})
		(a (e^{(a-b)t} s^K (1-s) + s^{K+1})-b s^K)}.
\end{eqnarray*}
where $G(s,t) - G(0,t)$ can be computed from~(\ref{eq:branch1}) using
elementary algebra.
The dominant term in the denominator of the above quantity is
$e^{(a-b)t} s^K (1-s)$, which achieves the maximum
\[ \frac{e^{(a-b)t}}{K+1} \left(1-\frac{1}{K+1}\right)^K
= \frac{e^{(a-b)t}}{K} \left(1-\frac{1}{K+1}\right)^{K+1} \]
at $s=K/(K+1)$. For sufficiently large $K$, this is at least
$e^{(a-b)t}/(3K)$. Therefore
\begin{eqnarray*}
	P(1 \le \xi(t) \le K) \le
	\frac{(a-b)^2 \frac K {K+1}}{(a - b e^{-(a-b)t}) \left(
		a \frac{e^{(a-b)t}}{3K} -b\right)}
	\le C_{a,b} \left(a \frac{e^{(a-b)t}}{3K}-b \right)^{-1},
\end{eqnarray*}
which implies the desired conclusion of (a.2), if $K \le e^{(a-b)t}/6$.

For (a.3), we observe that $M(t)=e^{-(a-b)t} \xi(t)$ is a martingale
with maximum jump size 1 and quadratic variation
$\langle M \rangle(t)=\int_0^t e^{-2(a-b)s} (a+b) \xi(s) \ ds$.
Burkholder's inequality implies
\begin{eqnarray*}
        E\left[ \sup_{s\le t} M(s) \right]
        &\le& C + C \int_0^t e^{-2(a-b)s} (a+b) E[\xi(s)] \ ds \\
        &=& C + C \int_0^t e^{-2(a-b)s} (a+b) e^{(a-b)s} \ ds
        \le C_{a,b}.
\end{eqnarray*}
Therefore
$E\left[ \sup_{s\le t} \xi(s) \right] \le C_{a,b} e^{(a-b)t}$,
which implies (a.3).

For (b), we observe that
\begin{eqnarray*}
        P(\xi^{(k)}(t)=0) = \left( 1-\frac{b-a}{b e^{(b-a)t}-a} \right)^k.
\end{eqnarray*}
For sufficiently large $t$, we have
\begin{eqnarray*}
	\frac{b e^{(b-a)t}-a}{b-a} = \frac{e^{(b-a)t}-\frac a b}{1-\frac a b}
	\ge e^{(b-a)t}-\frac a b
	\ge \frac 1 {1.1} e^{(b-a)t},
\end{eqnarray*}
therefore
\begin{eqnarray*}
        P(\xi^{(k)}(t)=0) &\ge& \left[ \left( 1-\frac{b-a}{b e^{(b-a)t}-a}
		\right)^{\frac{b e^{(b-a)t}-a}{b-a}} \right]^{
		k 1.1 e^{-(b-a)t}}
	\ge e^{-1.2 k e^{-(b-a)t}} \\
	&\ge& 1-1.2 k e^{-(b-a)t},
\end{eqnarray*}
if $t$ is sufficiently large and $k e^{-(b-a)t}$ is sufficiently small.

For (c), we observe that
$\xi^{(k)}=\xi^{(1)}_1+\xi^{(1)}_2+\ldots+\xi^{(1)}_k$, where
$\xi^{(1)}_i, i=1,\ldots,k$ are independent copies of $\xi^{(1)}$.
Therefore
\begin{eqnarray*}
        P(\xi^{(k)}(t)\in [1,K])
        \le P(\xi^{(1)}_i \in [1,K] \mbox{ for some } i=1,\ldots,k)
        \le k C_{a,b} K e^{-(a-b)t}
\end{eqnarray*}
by part (a.2) of this lemma. Part (d) is a direct consequence of part (c).
\qed

\bibliographystyle{plain}
\bibliography{ref2}

\end{document}